\documentclass[12pt]{article}

\usepackage{graphicx}
\newcommand{\beq}{\begin{equation}} 
\newcommand{\eeq}{\end{equation}}
\newcommand{\ba}{\begin{array}} 
\newcommand{\ea}{\end{array}}
\newcommand{\ben}{\begin{enumerate}}
\newcommand{\een}{\end{enumerate}} 
\newtheorem{theorem}{Theorem}
\newcommand{\bth}{\begin{theorem}} 
\newcommand{\eth}{\end{theorem}}
\newtheorem{definition}{Definition}
\newcommand{\bdf}{\begin{definition}}
\newcommand{\edf}{\end{definition}} 
\newtheorem{notation}{Notation}
\newcommand{\bnt}{\begin{notation}} 
\newcommand{\ent}{\end{notation}}
\newtheorem{proposition}{Proposition}
\newcommand{\bpr}{\begin{proposition}}
\newcommand{\epr}{\end{proposition}} 
\newtheorem{lemma}{Lemma}
\newcommand{\blm}{\begin{lemma}} 
\newcommand{\elm}{\end{lemma}}
\newtheorem{remark}{Remark}
\newcommand{\brm}{\begin{remark}} 
\newcommand{\erm}{\end{remark}}
\newtheorem{corollary}{Corollary}
\newcommand{\bcr}{\begin{corollary}} 
\newcommand{\ecr}{\end{corollary}}

\begin{document}
\title{Functional Analysis for Helmholtz Equation in the Framework of Domain Decomposition }

\author {Mikhael Balabane\footnote{ University Paris 13 - France -
  email: balabane@math.univ-paris13.fr}\date{Sept. 2008} }

\maketitle

\begin{abstract}\noindent
This paper\footnote{CMS 65M12, 35P10, 35P15, 35P25}  gives a geometric description of functional spaces  related to  Domain Decomposition techniques for computing solutions of Laplace and Helmholtz equations. Understanding the geometric structure of these spaces leads to  algorithms for solving  the equations. It leads also to a new interpretation of classical algorithms, enhancing convergence.  The algorithms are given and convergence is proved. 
\noindent This is done by building tools enabling geometric interpretations of the operators related to Domain Decomposition technique. The Despres operators, expressing conservation of energy for Helmholtz equation, are defined on the fictitious boundary and their spectral properties proved.It turns to be the key for proving convergence of the given algorithm for Helmholtz equation in a non-dissipating cavity.

\noindent Using these tools, one can prove that the Domain Decomposition setting for the Helmholtz equation leads to an ill-posed problem. Nevertheless, one can prove that if a solution exists, it is unique. And that the algorithm do converge to the solution.

\end{abstract}

\section {Introduction}

\noindent In the framework of domain decomposition,  given a bounded open set $\Omega=\Omega_1 \cup \Omega_2\cup \Gamma$ where the two  open sets $\Omega_1$ and $\Omega_2$ are not overlapping, and $\Gamma$  a common subset of their boundary (called the fictitious boundary), and given a (global) solution  $u$ of the Helmholtz equation 
$$\Delta u+ k^2u=f\in L^2(\Omega)\quad {\rm and}\quad u\in H^1_0(\Omega)
$$
the aim of this paper is  to understand the dynamics of the sequence $(v^n_1,v^n_2)_{n\in N}$ solving separatly the Helmholtz equations on $\Omega_1$ and $ \Omega_2$, when equating the fluxes through $\Gamma$: ($m=1,2$ resp. $m'=2,1)$
$$ {\partial v_m^n \over \partial n_m}-i\gamma v_m^n=-{\partial v_{m'}^{n-1}
\over \partial n_{m'}}-i\gamma v_{m'}^{n -1}\quad {\rm on}\quad \Gamma
$$
Its ultimate aim is to  prove convergence to $(u_{\vert \Omega_1},u_{\vert \Omega_2})$ of the sequence $(u^n_1,u^n_2)_{n\in N}$ solving the Helmholtz equations on $\Omega_1 $ and $ \Omega_2$ with a penalization on the boundary $\Gamma$ is added, namely:
$$ {\partial u_m^n \over \partial n_m}-i\gamma u_m^n=\theta [ {\partial u_m^{n-1} 
\over \partial n_m}-i\gamma u^2_{n-1}]-(1-\theta )[{\partial u_1^{n-1}
\over \partial n_{m'}}+i\gamma u_{m'}^{n -1}]\quad {\rm on}\quad \Gamma
$$

\noindent For this sake,  the geometry of the set of solutions of the Helmholtz equation on $\Omega_1 \times \Omega_2$ with equated energy fluxes is studied, through the study of the coupling  operator defined on $L^2(\Gamma) \times L^2(\Gamma)$ which intertwins the fluxes. It turns out that the key for understanding the convergence of the sequence   $(u^n_1,u^n_2)_{n\in N}$ is the analysis of the spectral properties of the intertwinnig operator.
\medskip

\noindent Using these tools, one can prove that the Domain Decomposition setting for the Helmholtz equation leads to an ill-posed problem. Nevertheless, one can prove that if a solution exists, it is unique. And that the algorithm do converge to the solution.
\medskip

\noindent Convergence of the penalized algorithm is proven and numerical tests for solving the Helmholtz equation through this domain decomposition algorithm are given.
\medskip

\noindent The geometric analysis given here provides the theoretical background for another numerical algorithm for computing the global solution $u$, by a specific spectral method. A forthcoming paper describes and gives the numerical analysis of this algorithm.
\medskip

 \noindent  This domain decomposition algorithm (in a   {\it dissipating cavity} case,  i.e. with a Sommerfeld-like radiation condition on part of the boundary), was first initiated and studied by B.Despres   in {\bf [D1]} {\bf [D2]} {\bf [BD]}, and  computational results given by J.D.Benamou  {\bf [B]} {\bf [BD]}, F.Collino and P.Joly {\bf [CGJ]}.
  \medskip

\noindent In order to perform the geometric analysis of the set of solutions of the Helmholtz equation on 
$\Omega_1 \times \Omega_2$, one has first to make a complete description of the geometry of the set of solutions of the Laplace equation on 
$\Omega_1 \times \Omega_2$. Geometric properties of this set proven below makes it possible to revisit the classical penalized Dirichlet/Neumann domain decomposition algorithm (with penalization) for solving the Laplace equation. A new version of this algorithm is given here, and proved to converge to the global solution, enhancing the usual assumption on the penalization parameter. 
  \medskip

\noindent This completes classical results  by O.Widlund  {\bf [PW]},  P.L.Lions {\bf [L]}, or A.Quar\-teroni and A.Valli {\bf [ FMQT]} {\bf [FQZ]} {\bf [QV]}. 
 \medskip

\noindent The paper is organized as follows: in  {\bf section 2} basic facts are revisited, although classical, and completed in order to set  the geometric framework needed. (It also makes the paper self contained). A precise study of duality, and the link with the Poincare-Steklov operators, is performed, which turns to be central for the remainder of the paper. In {\bf section 3} a new version of the Dirichlet/Neumann algorithm for the Laplace equation is given, and convergence is proved. In {\bf section 4} geometric tools for the Helmholtz equation, and  related domain decomposition algorithm, are given. Despres operators are studied and their spectral properties investigated. As is the intertwinnig operator. In {\bf section 5}, convergence of the domain decomposition algorithm for Helmholtz equation is proved. In  {\bf section 6}  numerical tests are given.

 \medskip

\noindent Throughout this paper, when dealing with  the Helmholtz equation, the frequency $k$  is assumed to be non-resonnant for the Dirichlet boundary condition. More precisely we shall always make the following

\noindent {\bf Assumption (A)} $-k^2$ is not an eigenvalue of the Laplace operator on $\Omega$ with Dirichlet boundary condition, i.e. the following problem is well posed for $f\in L^2(\Omega)$:
$$\Delta u+ k^2u=f\quad {\rm and}\quad u\in H^1_0(\Omega) $$

\noindent We shall also adopt the following

\noindent {\bf Notation (N)} normal derivatives at the boundary of an open set are always meant as the derivative along the {\bf outward} unit normal vector

\section {Basics}
\medskip
Let $ \Omega  \subset R^d$ be a bounded open set   whose boundary  $\partial \Omega$ is a $ C^1 $-submanifold of $R^d$. Let $\Gamma$ be an open  $ C^\infty $-submanifold of $R^d$,  such that:
$$\Omega=\Omega_1\cup \Omega_2\cup \Gamma,\,\,\, \partial \Omega_1=(\partial \Omega \cap \partial \Omega_1)\cup \Gamma,\,\,\,\partial \Omega_2=(\partial \Omega \cap \partial \Omega_2)\cup \Gamma$$  
where $\Omega_1$ and $\Omega_2$ are open sets in $R^d$.  We assume  that $\Omega_1$ and $\Omega_2$ fulfill the strict cone property (see {\bf [Ag]} for instance) and that $\Gamma$ is transverse to $\partial \Omega$ in the following sense:
 $\overline \Gamma$ is a $ C^1 $-submanifold of $R^d$ with boundary, and there exists $a< 1$  such that for any $\sigma \in \partial \Omega \cap {\overline \Gamma}$, we have: 
\beq \label {XF1}  -a\leq n_{\overline \Gamma}(\sigma).n_{\partial \Omega}(\sigma)\leq a
\eeq
where $n_{\overline \Gamma}(\sigma) \in C^0(\overline \Gamma)$ is a  unit vector normal to $\overline \Gamma$ at $\sigma$ and $n_{\partial \Omega}(\sigma) \in C^0(\partial \Omega)$ a  unit vector normal to ${\partial \Omega}$ at $\sigma$.

\subsection{Functional spaces associated to $\Gamma$}
Let $H_0^1(\Omega)$ be endowed with the scalar product
$$(u,v)_{H_0^1(\Omega)}=\int_\Omega\nabla u\nabla {\overline v}dx$$
For $m=1,2$ let $H_m=\{u\in H^1(\Omega_m);u_{\vert \partial \Omega \cap \partial \Omega_m}=0\}$.
\noindent Boundedness of the trace operators from $H^1_0(\Omega)$ to  $H^{1\over 2}_0(\partial \Omega \cap \partial \Omega_m)$ imply that these are Hilbert spaces when endowed with the scalar products:
$$(u,v)_{H_m}=\int_{\Omega_m}\nabla u\nabla {\overline v}dx$$
Let $\rho^\Gamma$ (resp. $\rho_m^\Gamma$ for $m=1,2$) be the trace operator on $\Gamma$, i.e. the bounded linear operator from $H_0^1(\Omega)$ (resp. $H_m$) to $H^{1/2}(\Gamma)$ which maps $u$ to $u_{\vert \Gamma}$. 

 \noindent Let $$\Lambda=\{u_{\vert \Gamma}; u\in H_0^1(\Omega)\}=H_0^1(\Omega)/Ker \rho^\Gamma \simeq (Ker \rho^\Gamma)^\bot$$
 and for $m=1,2$
 $$  \Lambda_m=\{u_{\vert \Gamma}; u\in H_m\}=H_m/Ker \rho_m ^\Gamma\simeq (Ker \rho_m^\Gamma)^\bot$$
\brm \label{X0} Obviously $Ker \rho_m^\Gamma=H^1_0(\Omega_m)$, $ \Lambda_m\subset  H^{1/2}(\Gamma)$, $  \Lambda \subset H^{1/2}(\Gamma)$
\erm
\medskip

\noindent Because $\rho^\Gamma$ and $\rho_m^\Gamma$ are  bounded, $\Lambda$ and $\Lambda_m$ are Hilbert spaces when endowed with the following norms:
$$\forall \lambda \in \Lambda,\,\,\,\Vert \lambda \Vert_\Lambda=\inf_{\{u;u_{\vert \Gamma}=\lambda\}}\Vert u\Vert_{H_0^1(\Omega)}\quad {\rm and}\quad \forall \lambda \in \Lambda_m,\,\,\,\Vert \lambda \Vert_{\Lambda_m}=\inf_{\{u;u_{\vert \Gamma}=\lambda\}}\Vert u\Vert_{H_m}$$

\medskip

\brm \label{X00} Obviously, for any $ v\in H_0^1(\Omega)$ and $w\in H_m$
$$ \Vert \rho^\Gamma(v)\Vert_\Lambda\leq  \Vert v\Vert_{H^1_0(\Omega)} \quad and \quad \Vert \rho_m^\Gamma(w)\Vert_\Lambda\leq  \Vert w\Vert_{H_m}$$
\erm

\medskip

\bpr \label{X1} Let $ m=1,2$.   For any $\lambda \in \Lambda_m, \mu \in \Lambda_m$:

1- There exists a unique $u^\lambda_m \in H_m$ such that $\Delta  u^\lambda_m=0$ in $\Omega_m$ and $\rho_m^\Gamma(u_m^\lambda)=\lambda$

2- One has: $\Vert \lambda \Vert_{\Lambda_m}=\Vert u^\lambda_m \Vert_{H_m}$ and $(\lambda,\mu)_{\Lambda_m}=(u_m^\lambda,u_m^\mu)_{H_m}$

3- For any $ \lambda \in \Lambda$ let $u^\lambda=u^\lambda_m$ on $\Omega_m$,  $m=1,2$. Then
$$  \Vert \lambda \Vert_{\Lambda}^2=\Vert u^\lambda_1 \Vert_{H_1}^2+\Vert u^\lambda_2 \Vert_{H_2}^2 $$

4- Using the previous notation, for any  $ \lambda \in \Lambda$, $ \mu \in \Lambda$
$$\quad (\lambda,\mu)_{\Lambda}=(u^\lambda,u^\mu)_{H^1_0(\Omega)}$$
\epr
\noindent {\it proof:}

1- Uniqueness follows well posedness of the Laplace problem in $H_0^1(\Omega_m)$. In order to prove existence, let $u\in H_m$ be such that $ \lambda=\rho_m^\Gamma(u)$. Then $\Delta u \in H^{-{1}}(\Omega_m)$.  Let $v$ be the unique solution in  $H_0^1(\Omega_m)$ of $\Delta v=\Delta u$. Then $u_m^\lambda=u-v$ fulfills the property.

2- Because of  remark \ref{X0}, one has to show:
$$w\in H^1_0(\Omega_m)\Rightarrow (u_m^\lambda,w)_{H_m}=0$$
which follows from the Green formula. 

3- Because $\rho_1^\Gamma(u^\lambda_1)=\rho_2^\Gamma(u^\lambda_2)=\lambda$ one has  $u^\lambda \in H_0^1(\Omega)$. Obviously it is orthogonal to $Ker \rho^\Gamma$. so 
$$ \Vert \lambda \Vert_{\Lambda}^2=\Vert u^\lambda \Vert_{H^1(\Omega)}^2=\Vert u^\lambda_1 \Vert_{H_1}^2+\Vert u^\lambda_2 \Vert_{H_2}^2$$

\medskip
\subsection{$\Lambda=\Lambda_1=\Lambda_2$}

By symmetry, it is enough to prove  $\Lambda=\Lambda_1$. In order to prove this algebraic and topological  equality,  two key tools are needed. The first tool is the Calderon extension theorem  {\bf [Ag]},  which applies here because $\Omega_1$ has  the strict cone property, by assumption, and which gives a bounded linear operator $E$ from $H^1(\Omega_1)$ to $H^1(R^d)$ such that:
$$\forall w\in H^1(\Omega_1),\quad E w_{\vert \Omega_1}=w$$
The second key tool is:
\bth \label{X2}There exists a bounded linear operator $\tau$ in $H^1(R^d)$ such that:
$$\forall v\in H^1(R^d), \quad v_{\vert \partial \Omega_1 \cap \partial \Omega}=0\Rightarrow (\tau v_{\vert  \partial \Omega}=0\quad {\rm and}\quad \tau v_{\vert  \Omega_1}=v_{\vert  \Omega_1})$$
\eth
\noindent {\it proof:}

\noindent Assumption (\ref{XF1}) gives a finite open covering $(\omega_j)_j$ of $\overline \Gamma \cap \partial \Omega$ and a change of variables $(a^j)_j$ such that $V_j=a^j(\omega_j)$ is a neighbourhood of zero in $R^d$ and:   $$a^j(\Gamma \cap \omega_j)=\{z^j\in V_j;z^j_1=0,z^j_2>0\}\quad  and\quad a^j(\partial \Omega \cap \omega_j)=\{z^j\in V_j;z^j_2=0\}$$
$$a^j( \Omega_1 \cap \omega_j)=\{z^j\in V_j;z^j_1<0,z^j_2>0\}$$
$$ a^j(\Omega_2 \cap \omega_j)=\{z^j\in V_j;z^j_1>0,z^j_2>0\}
$$
Regularity of the submanifold $\Gamma$ gives an open covering $(\omega'_k)_k$ of $\Gamma$ and change of variables $(b^k)_k$ such that $W_k=b^k(\omega'_k)$ is a neighbourhood of zero in $R^d$ and:
$$b^k(\Gamma \cap \omega'_k)=\{z^k\in W_k;z^k_1=0\}$$
$$ b^k(\Omega_1\cap \omega'_k)=\{z^k\in W_k;z^k_1<0\}\quad and \quad b^k(\Omega_2\cap \omega'_k)=\{z^k\in W_k;z^k_1>0\}
$$   
Compactness of $\overline \Gamma$ enables to select a finite subcovering of $\overline \Gamma$ still denoted by $(\omega_j)_j\cup (\omega'_k)_k$ having the previous properties.
\medskip

\noindent Let $\omega^1=R^d\setminus \overline\Omega_1$ and $\omega^2=R^d\setminus \overline\Omega_2$, so:
$$R^d= \omega^1\cup \omega^2\cup (\cup_j  \omega_j)\cup (\cup_k  \omega'_k)$$
Let $$(\alpha^1,\alpha^2,(\alpha_j)_j,(\alpha'_k)_k)$$ be a $C^\infty$-partition of unity associated with this open covering of $R^d$.

\noindent Let $\varepsilon>0$ be such that $ \forall z\in (\cup_jV_j)\cup (\cup_k W_k) \quad 0<z_1^j<\varepsilon\Rightarrow z\in\Omega_2$

\noindent Let  $\psi(s) \in C^\infty(R)$ be equal to one for $s<0$ and zero for $s>\varepsilon$. 

\noindent Let $\varphi(\theta) \in C^\infty(R)$ be equal to zero for $\theta<0$ and equal to one for $\theta>{\pi\over 2}$
\medskip

\noindent  For any $v\in H^1(R^d)$,
$$v=\alpha^1v+\alpha^2v+\sum_j\alpha_jv+\sum_k\alpha'_kv
$$

\noindent we define $\tau v$ as:
$$\tau(v)=\tau(\alpha^1v)+\tau(\alpha^2v)+\sum_j\tau(\alpha_jv)+\sum_k\tau(\alpha'_kv)$$
with:

\noindent  $\tau(\alpha^1v)\equiv 0$

\noindent  $\tau(\alpha^2v)= \alpha^2v$

\noindent  $\tau(\alpha'_kv)(z^k)=\psi(z^k_1)\alpha'_k(z^k)v(z^k)$ in the local coordinates.

\noindent These three quantities are multiplication of $v$ by $C^\infty$ functions, which are bounded  as well as all their derivatives. It is linear and bounded in $H^1(R^d)$ with respect to $v\in H^1(R^d)$. 

\noindent In order to define $\tau(\alpha_jv)$, we first write $\alpha_jv$ in the cylindrical coordinates as follows:
$$\alpha_jv(z^j_1,z^j_2,z^j_3,..,z^j_d)=\widetilde{\alpha_jv}(r^j,\theta^j,z^j_3,..,z^j_d)\quad {\rm with}\quad z^j_1=r^jcos(\theta^j),\,z^j_2=r^jsin(\theta^j)
$$
and define $\tau(\alpha_jv)$ in these coordinates as:
$$\widetilde{\tau(\alpha_jv)}(r^j,\theta^j,z^j_3,..,z^j_d)=\varphi(\theta^j)\widetilde{\alpha_jv}(r^j,\theta^j,z^j_3,..,z^j_d)$$
This quantity is linear with respect to $v$, and we prove its boundedness in $H^1(R^d)$ with respect to $v\in H^1(R^d)$ as follows (we omit the index $j$ and denote the measure $dz_3...dz_d$ by $d\overline{z}$):
$$\Vert \tau(\alpha v)\Vert_{H^1(R^d)}^2=\int_V\vert \widetilde{\tau(\alpha v)}\vert^2rdrd\theta d\overline{z}+\int_V\vert {\partial \over \partial r }\widetilde{\tau(\alpha v)}\vert^2rdrd\theta d\overline{z}$$
$$+\int_V {1\over r^2}\vert {\partial \over \partial \theta } \widetilde{\tau( \alpha v)}\vert^2rdrd\theta d\overline{z}+\sum_3^d\int_V\vert {\partial \over \partial z_j }\widetilde{\tau(\alpha v)}\vert^2rdrd \theta d\overline{z}
$$
$$=\int_V\vert \varphi(\theta)\widetilde{\alpha v}\vert^2rdrd\theta d\overline{z}+\int_V\vert  \varphi(\theta){\partial \over \partial r }\widetilde{\alpha v}\vert^2rdrd\theta d\overline{z}$$
$$+\int_V {1\over r^2}\vert  \varphi(\theta){\partial \over \partial \theta } \widetilde{\alpha v}+ \varphi'(\theta) \widetilde{\alpha v}\vert^2rdrd\theta d\overline{z}+\sum_3^d\int_V\vert  \varphi(\theta) {\partial \over \partial z_j }\widetilde{\alpha v}\vert^2rdrd\theta d\overline{z}
$$
$$\leq sup\vert \varphi\vert^2[\int_V\vert \widetilde{\alpha v}\vert^2rdrd\theta d\overline{z}+\int_V\vert  {\partial \over \partial r }\widetilde{\alpha v}\vert^2rdrd\theta d\overline{z}$$
$$+\int_V {2\over r^2}\vert {\partial \over \partial \theta } \widetilde{\alpha v}\vert^2rdrd\theta d\overline{z}+\sum_3^d\int_V\vert   {\partial \over \partial z_j }\widetilde{\alpha v}\vert^2rdrd\theta d\overline{z}]$$
$$+2sup\vert \varphi'\vert^2\int_V {1\over r^2}\vert  \widetilde{\alpha v}\vert^2rdrd\theta d\overline{z}
$$
$$\leq 2sup\vert \varphi\vert^2\Vert \alpha v\Vert^2_{H^1(R^d)}+2sup\vert \varphi'\vert^2\int_V {1\over r^2}\vert  \widetilde{\alpha v}\vert^2rdrd\theta d\overline{z}
$$
In order to estimate this last  quantity 
we use the assumption $v_{\vert \partial \Omega_1 \cup \partial \Omega}=0$ to have:
$$ \widetilde{\alpha v}(r,\theta,\overline{z})=-\int_\theta^\pi{\partial \over \partial \theta }\widetilde{\alpha v}(r,s,\overline{z})ds$$
which gives (with $\varepsilon '$ the radius of the support of $\alpha$ in the $r$ variable, and $B$ a ball containing the support of $\alpha$ in the $\overline{z}$ variable):

$$\int_V {1\over r^2}\vert  \widetilde{\alpha v}\vert^2rdrd\theta d\overline{z}
\leq \int_B \int _{0}^{\varepsilon'}\int_{-\pi}^\pi   {1\over r^2}\vert \int_\theta^\pi {\partial \over \partial \theta}\widetilde{\alpha v}(r,s,\overline{z})ds\vert^2rdrd\theta d\overline{z}$$
$$\leq 4\pi^2 \int_B \int _{0}^{\varepsilon'}\int_{-\pi}^\pi   {1\over r^2}\vert  {\partial \over \partial \theta}\widetilde{\alpha v}(r,s,\overline{z})\vert^2rdrdsd\overline{z}\leq 4\pi^2\Vert \alpha v\Vert^2_{H^1(R^d)}
$$
we summarize to have:
$$\Vert \tau(\alpha v)\Vert^2_{H^1(R^d)}\leq C\Vert \alpha v\Vert^2_{H^1(R^d)}\leq C'\Vert  v\Vert^2_{H^1(R^d)}
$$
and this ends the proof of the   boundedness of $\tau$ in $H^1(R^d)$.

\medskip
\noindent We end the proof of theorem \ref{X2} using the following obvious observations:

$\tau v_{\vert \Omega_1}= v_{\vert \Omega_1}$ because $\psi\equiv 1$ on $R_-$ and $\varphi\equiv 1$ for $\theta\geq {\pi \over 2}$

$\tau v_{\vert \partial \Omega_1\cup \partial \Omega}= 0$ by assumption

$\tau v_{\vert \partial \Omega_2\cup \partial \Omega}= 0$ because $\varphi(0)=0$ and $\psi\equiv 0$ for  $z_1^j>\varepsilon$.

\bcr \label{X3} For $m=1,2$, there exists a bounded linear map $E_m$ from $H_m$ to $H^1_0(\Omega)$ such that
$$\forall u \in H_m\quad (E_m{u})_{\vert \Omega_m}\equiv u$$
\ecr
\noindent {\it proof:}  for $m=1$ for instance  let $E_1u$ be the restriction to $\Omega$ of $\tau E u$.Boundedness follows from theorem \ref{X2}.

\bcr \label{X4}  $\Lambda_1=\Lambda_2=\Lambda$ and the three norms $\Vert .\Vert_{\Lambda},\Vert .\Vert_{\Lambda_1}$ and $\Vert .\Vert_{\Lambda_2}$ are equivalent.
\ecr
{\it proof:} obviously $\Lambda \subset \Lambda_m$ and the previous corollary gives the converse inclusion. Moreover we have $\Vert .\Vert_{\Lambda}\geq\Vert .\Vert_{\Lambda_m}$ and  the previous corollary gives the converse inequality.

\bcr \label{X5}
${\cal D}(\Gamma)$ is a dense subspace of $\Lambda$ for any of the three norms.
\ecr
{\it proof:}  ${\cal D}(\Gamma)$ is the set of traces on $\Gamma$ of  functions in ${\cal D}(\Omega)$ because $\Gamma$ is a $C^\infty$submanifold. Density of ${\cal D}(\Gamma)$ for the $\Vert .\Vert_{\Lambda}$ norm follows density of ${\cal D}(\Omega)$ in $H^1_0(\Omega)$. Equivalence of the three norms ends the proof.

\medskip
\brm \label{X000} If we denote as usual by $H^{1\over 2}_0(\Gamma)$ the closure of  ${\cal D}(\Omega)$ in $H^{1\over 2}(\Gamma)$ (which exists because $\Gamma$ is $C^\infty$),  then the previous corollary asserts that $\Lambda\subset H^{1\over 2}_0(\Gamma)$. Boundedness of the trace operators gives constants $C,C_m$ such that:
$$\forall \lambda \in \Lambda,\,\,\Vert \lambda \Vert_{H^{1\over 2}(\Gamma)}\leq C\Vert \lambda \Vert_\Lambda \quad {\rm and}\quad \forall \lambda \in \Lambda_m,\,\,\,,\Vert \lambda \Vert_{H^{1\over 2}(\Gamma)}\leq C_m\Vert \lambda \Vert_{\Lambda_m}$$
\erm

\subsection{Well-posedness of the Laplace-Dirichlet problem in $H^{-1}(\Omega_m)\times\Lambda$}

\bth \label{X6} For $m=1,2$

\noindent 1- for any $(f,\lambda)\in H^{-1}(\Omega_m)\times\Lambda$ there exists a unique $u\in H_1$ such that $$\Delta u=f \,\, in\,\,  \Omega_m\quad and \quad \rho_m^\Gamma(u)=\lambda$$

\noindent 2-  we have the estimate $\Vert u\Vert_{H_m}\leq \Vert f\Vert_{H^{-1}(\Omega_m)}+\Vert \lambda\Vert_{\Lambda_m}$
\eth
{\it proof:} Let $u_m^\lambda$ be given by Proposition \ref{X1}. Let $v=u-u_m^\lambda$. The problem is equivalent to
$$v\in H^1_0(\Omega_m)\quad and \quad \Delta v=f$$
This is a well-posed problem and we have, because the Riesz representation operator is isometric,
$\Vert v\Vert_{H^1_0(\Omega_m)}=\Vert f\Vert_{H^{-1}(\Omega_m)}$. So
$$\Vert u\Vert_{H_m}\leq\Vert u_m^\lambda\Vert_{H_m}+\Vert v\Vert_{H_m}\leq\Vert f\Vert_{H^{-1}(\Omega_m)}+\Vert \lambda\Vert_{\Lambda_m}$$

\subsection{Duality and the Poincare-Steklov operators}
Let $\Lambda'$ denote the dual space to $\Lambda$, endowed with one of the three equivalent norms associated with the equivalent norms on $\Lambda$ defined previously.

\noindent We denote by $(.,.)_{\Lambda \Lambda'}$ the duality product, and by  $(.,.)_{{\cal D} {\cal D'}}$ the duality product in ${\cal D'}(\Gamma)$. 

\noindent Because of Corollary \ref{X5} and remark \ref{X000}, we have the usual injections:
$${\cal D} \subset \Lambda \subset H^{1\over 2}_0(\Gamma) \subset L^2(\Gamma)\subset H^{-{1\over 2}}(\Gamma)\subset \Lambda' \subset {\cal D'}
$$
and for any $\lambda \in {\cal D}(\Gamma)$ and $\nu \in \Lambda' $: $(\nu,\lambda)_{\Lambda \Lambda'}=(\nu,\lambda)_{{\cal D} {\cal D'}}=(\nu,\overline{\lambda})_{L^2(\Gamma)}$

\noindent where  $\overline{\eta}$ denotes the complex conjugate of functions or distributions $\eta$. 

\noindent For any of the three scalar products on $\Lambda$ we have  $(\lambda,\eta)=\overline{(\overline{\lambda},\overline{\eta})}$ so for all norms:
$$\forall \lambda \in \Lambda, \quad\Vert \lambda \Vert=\Vert \overline{\lambda}\Vert \quad and \quad  \forall \nu \in \Lambda',\quad  \Vert \nu \Vert=\Vert \overline{\nu}\Vert $$

\bnt \label{X66} Let $\tilde S$ denote the antilinear Riesz representation operator for $\Lambda'$ in the $\Lambda$-scalar product, and  $\tilde S_m$ ($m=1,2$) this representation in the $\Lambda_m$-scalar product, i.e.
$$\forall \lambda \in \Lambda,\forall \nu \in \Lambda'\quad(\nu,\lambda)_{\Lambda \Lambda'}=( \lambda,\tilde S^{-1}\nu)_{\Lambda}=( \lambda,\tilde S^{-1}_1\nu)_{\Lambda_1}=( \lambda,\tilde S^{-1}_2\nu)_{\Lambda_2}
$$
Let  $S$ (resp. $S_m$) be the linear isometric bijections from $\Lambda$ to $\Lambda'$ defined by
$$\forall \Lambda \in \Lambda \quad \tilde S \lambda=S \overline{\lambda},\quad  \tilde S_1 \lambda=S_1 \overline{\lambda}, \quad \tilde S_2 \lambda=S_2 \overline{\lambda}\quad$$
\ent
\noindent  We denote by  $n_m $ the normal unit  vector on $\Gamma$  pointing outward with respect to $\Omega_m$. We denote by ${\partial \varphi \over \partial n_m}$ the normal derivative on $\Gamma$ of $\varphi \in {\cal D}(\Omega)$, and by ${\partial  \over \partial n_m}$ bounded extensions of this operator to any functional space. 

\medskip

\noindent We use  Proposition \ref{X1} to have:

\bpr \label{X7} $$\forall \lambda \in \Lambda \quad S_m\lambda={\partial u_m^\lambda\over \partial n_m} \quad and \quad S=S_1+S_2$$
\epr
{\it proof:} by Green formula
$$\forall \lambda \in \Lambda,\eta \in \Lambda\quad  (\tilde{S_m} {\eta},\lambda)_{\Lambda_m\Lambda'_m}= (\lambda,\eta)_{\Lambda_m}=$$
$$(u_m^\lambda,u_m^\eta)_{H^1(\Omega_m)}=\int_{\Omega_m}\nabla u_m^\lambda \overline {\nabla u_m^\eta}dx=\int_\Gamma\lambda 
 \overline {\partial u_m^\eta \over \partial n_m}d\sigma$$
 this proves that the distribution ${\partial u_m^\eta \over \partial n_m}$ is bounded in the $\Lambda_m$ norm, so   ${\partial u_m^\eta \over \partial n_m} \in \Lambda'$ and 
 $$\tilde{S_m} {\eta}=  \overline {\partial u_m^\eta \over \partial n_m}=  \overline {\partial{u_m^\eta }\over \partial n_m}={\partial u_m^ {\overline {\eta }}\over \partial n_m}
$$
\bcr  \label{X8}$$\overline{S \lambda}=S\overline{ \lambda}\quad \overline{S_m \lambda}=S_m\overline{ \lambda}
$$
\ecr

\brm \label{X29} For any $v\in H_m$ such that $\Delta v \in L^2(\Omega_m)$, we have ${\partial v\over \partial n_m}\in \Lambda'$ and 
$$\Vert {\partial v\over \partial n_m}\Vert_{\Lambda'}\leq C(\Vert v\Vert_{H_m }+\Vert \Delta v\Vert_{L^2(\Omega_m)})
$$
\erm
\noindent This is because for any $\varphi \in {\cal D}(\Omega)$ we have:
$$({\partial v\over \partial n_m},\varphi_{\vert \Gamma})_{{\cal D}{\cal D'}}=\int_{\Omega_m}\nabla \varphi \nabla vdx+\int_{\Omega_m} \varphi \Delta vdx
$$
\noindent and this formula shows that the distribution ${\partial v\over \partial n_m}$ is bounded on  $\Lambda$.

\subsection{Adjoints}
\bnt \label{X9}:

\noindent 1- For a bounded linear operator $T$ from $\Lambda$ to $\Lambda'$, we denote by $T'$ its adjoint for the $(\Lambda,\Lambda')$ duality, i.e.
$$\forall \lambda \in \Lambda, \forall \eta \in \Lambda\quad (T\eta, \lambda)_{\Lambda\Lambda'}= (T'\lambda, \eta)_{\Lambda\Lambda'}$$

\noindent 2- For a bounded linear operator $T$ from $\Lambda$ to $\Lambda$, we denote by $T^*$ the adjoint operator in $\Lambda$, i.e.
$$\forall \lambda \in \Lambda, \forall \eta \in \Lambda\quad (T\eta, \lambda)_{\Lambda}= ( \eta,T^*\lambda)_{\Lambda}$$

\ent  

\bpr
\label{X10}For $m=1,2$
$$S'=S\quad S'_m=S_m$$
\epr
\noindent {\it proof:} by definition of $S$ we have: $\forall \lambda \in \Lambda, \forall \eta \in \Lambda$

 $(S\lambda,\eta)_{\Lambda\Lambda'}=(\tilde {S}\overline{\lambda},\eta)_{\Lambda\Lambda'}=(\eta, \overline{\lambda})_\Lambda= \overline{( \overline{\lambda},\eta)}_\Lambda=\overline{( S\overline{\eta}, \overline{\lambda})}_{\Lambda,\Lambda'}={( S{\eta}, {\lambda})}_{\Lambda,\Lambda'}$

\bth
\label{X11}For $m=1,2$ let $m'=2,1$. For all $\lambda \in \Lambda, \eta \in \Lambda$

1- $(S_m^{-1}S_{m'}\eta,\lambda)_{\Lambda_m}=(\eta,\lambda)_{\Lambda_{m'}}$

2- $(S_m^{-1}S_{m'}\eta,\lambda)_{\Lambda_{m'}}=(S_{m'}\eta,S_{m'}\lambda)_{\Lambda'_{m}}$

3- $S_m^{-1}S_{m'}$ is self adjoint in $\Lambda_m$ and in $\Lambda_{m'}$

4- $(S_1^{-1}S_2+S_2^{-1}S_1)$ is selfadjoint in $\Lambda$

\eth
{\it proof:} 

\noindent 1- $(S_m^{-1}S_{m'}\eta,\lambda)_{\Lambda_m}=\overline{(\lambda,S_m^{-1}S_{m'}\eta)_{\Lambda_m}}=\overline{(S_{m'}\overline{\eta},\lambda)_{\Lambda\Lambda'}}=\overline{(\lambda,\eta)_{\Lambda_{m'}}}=(\eta,\lambda)_{\Lambda_{m'}}
$

\noindent 2- We use  Proposition \ref{X11} to have

$(S_m^{-1}S_{m'}\eta,\lambda)_{\Lambda_{m'}}=(S_{m'}\overline{\lambda},S_m^{-1}S_{m'}\eta)_{\Lambda\Lambda'}=(S_{m'}\eta,S_m^{-1}S_{m'}\overline{\lambda})_{\Lambda\Lambda'}=$

$(S_mS_m^{-1}S_{m'}\eta,S_m^{-1}S_{m'}\overline{\lambda})_{\Lambda\Lambda'}= (S_m^{-1}S_{m'}\overline{\lambda},S_m^{-1}S_{m'}\overline{\eta})_{\Lambda_m}=(S_{m'}\eta,S_{m'}\lambda)_{\Lambda'_{m}}$

\noindent 3- We use  Proposition \ref{X11} to have

$(S_m^{-1}S_{m'}\eta,\lambda)_{\Lambda_m}=\overline{(S_{m'}\overline{\eta},\lambda)_{\Lambda\Lambda'}}=\overline{(S_{m'}\lambda,\overline{\eta})_{\Lambda\Lambda'}}=\overline{(\overline{\eta},S_m^{-1}S_{m'}\overline{\lambda})_{\Lambda_m}}={({\eta},S_m^{-1}S_{m'}{\lambda})_{\Lambda_m}}$

On the other hand

$(S_m^{-1}S_{m'}\eta,\lambda)_{\Lambda_{m'}}=(S_{m'}\overline{\lambda},S_m^{-1}S_{m'}\eta)_{\Lambda\Lambda'}=(S_{m'}\eta,S_m^{-1}S_{m'}\overline{\lambda})_{\Lambda\Lambda'}=(S_m^{-1}S_{m'}\overline{\lambda},\overline{\eta})_{\Lambda_{m'}}=$ 

$(\eta,S_{m}^{-1}S_{m'}\lambda)_{\Lambda_{m'}}$

\noindent 4- follows 3

\bcr
\label{X12}{\it Coerciveness:} For $m=1,2$ let $m'=2,1$. There exists $C>0$ such that for all $\lambda \in \Lambda$:

1- $(S_m^{-1}S_{m'}\lambda,\lambda)_{\Lambda_m}=\Vert \lambda \Vert_{\Lambda_{m'}}^2\geq C \Vert \lambda \Vert_{\Lambda_{m}}^2$

2- $(S_m^{-1}S_{m'}\lambda,\lambda)_{\Lambda_{m'}}=\Vert S_{m'}\lambda \Vert_{\Lambda'_{m}}^2\geq C \Vert \lambda \Vert_{\Lambda_{m'}}^2$

3- $((S_1^{-1}S_2+S_2^{-1}S_1)\lambda,\lambda)_{\Lambda}\geq (1+C) \Vert \lambda \Vert_{\Lambda}^2$
\ecr

\subsection{On the Neumann problem}
\bpr
\label{X15} For $m=1,2$ and for any $\nu \in \Lambda'$ there exists a unique $u_m \in H_m$ such that 
$$\Delta u_m=0\quad in \quad \Omega_m\quad and \quad {\partial u_m\over \partial n_m}=\nu \quad on \quad \Gamma$$
Moreover
$$\Vert u_m\Vert_{H_m}=\Vert \nu\Vert_{\Lambda'}
$$
\epr
{\it proof:} uniqueness is straightforward, and existence is provided by proposition \ref{X1} and $u_m=u^{S_m^{-1}\nu}$

\noindent Moreover isometry of the Riesz representation gives:
$$\Vert  u_m\Vert_{H_m}=\Vert S_m^{-1}\nu\Vert_{\Lambda}=\Vert \nu\Vert_{\Lambda'}$$

\bpr \label{X16} For $m=1,2$ and for any $f\in L^2(\Omega_m)$ there exists a unique $u_m \in H_m$ such that 
$$\Delta u_m=f\quad {\rm and}\quad {\partial u_m\over \partial n_m}=0\quad {\rm on}\,\,\Gamma $$
Moreover
$$\Vert u_m\Vert_{H_m}\leq C\Vert f\Vert_{L^2(\Omega_m)}
$$
\epr
{\it proof:} Uniqueness is straightforward. For existence let $v_m\in H^1_0(\Omega_m)$ be the unique solution of $\Delta v_m=f$. In remark \ref{X29} we have shown that  $\nu={\partial v_m\over \partial n_m}\in \Lambda'$ and 
$$\Vert \nu\Vert_{\Lambda'_m}\leq C(\Vert f\Vert_{L^2(\Omega_m)}+ \Vert v_m\Vert_{H_m})\leq C(\Vert f\Vert_{L^2(\Omega_m)}+ \Vert f\Vert_{H^{-1}(\Omega_m)})\leq C\Vert f\Vert_{L^2(\Omega_m)}
$$
The function $u_m=v_m-u^{S_m^{-1}\nu}$ solves the problem, and we have:
$$\Vert u_m\Vert_{H_m}\leq \Vert v_m\Vert_{H_m}+\Vert u^{S_m^{-1}\nu}\Vert_{H_m}\leq
 \Vert f\Vert_{L^2(\Omega_m)}+\Vert S_m^{-1}\nu\Vert_{\Lambda_m}
 $$
 $$
\leq  \Vert f\Vert_{L^2(\Omega_m)}+\Vert \nu\Vert_{\Lambda'_m}
 \leq C\Vert f\Vert_{L^2(\Omega_m)}
$$

\section{A two-sided Dirichlet-Neumann domain decomposition algorithm for the Laplace operator}
\bpr
\label{X13} Let $f \in L^2(\Omega)$ and $u\in H^1_0(\Omega)$ be the unique solution of $ \Delta u=f $.  For $m=1,2$ let $f_m=f_{\vert \Omega_m}$ and $g_m\in H^1_0(\Omega_m)$ be the unique solution of  $ \Delta g_m=f_m$.  Let $\eta_m={g_m}_{\vert \Gamma}$
$$\lambda=u_{\vert \Gamma}\Longleftrightarrow (S_1+S_2)\lambda=-(S_1\eta_1+S_2\eta_2)$$
\epr

\noindent {\it proof:} the direct implication is stating continuity of $u$ and its normal derivatives through $\Gamma$. The converse implication states that taking $u_{\vert \Omega_m}=u_m^{\lambda}$ solves the global problem.

\medskip
\noindent We use the same notation as in the previous proposition to state:
\bth
\label{X14} Let $0<\theta<1$ with $ \Vert S_1^{-1}S_2+S_2^{-1}S_1\Vert_{{\cal L}(\Lambda)}<{2(1-\theta)\over \theta}$. Any sequence $(\lambda_n)_n\subset \Lambda$ which fulfills
  $$\lambda_{n+1}=((1-\theta)Id-{\theta\over 2} (S_1^{-1}S_2+S_2^{-1}S_1))\lambda_n-{\theta \over 2}(S_1^{-1}+S_2^{-1}))(S_1\eta_1+S_2\eta_2)$$
 do converge in $\Lambda$ (with geometric rate $1-{3\theta \over 2}$ at least) and its limit is $u_{\vert \Gamma}$.
\eth
{\it proof:} theorem \ref{X11} states selfadjointness of $S_1^{-1}S_2+S_2^{-1}S_1$ in $\Lambda$ and theorem \ref{X12} states coerciveness of $S_1^{-1}S_2+S_2^{-1}S_1$ so:
$$\Vert (1-\theta)Id-{\theta\over 2} (S_1^{-1}S_2+S_2^{-1}S_1)\Vert_{{\cal L}(\Lambda)}=$$
$$\sup_\Lambda{\vert (((1-\theta)Id-{\theta\over 2} (S_1^{-1}S_2+S_2^{-1}S_1))\lambda,\lambda)_\Lambda\vert \over \Vert \lambda\Vert^2_\Lambda}=
$$
$$\sup_\Lambda{ (((1-\theta)Id-{\theta\over 2} (S_1^{-1}S_2+S_2^{-1}S_1))\lambda,\lambda)_\Lambda \over \Vert \lambda\Vert^2_\Lambda}\leq 1-\theta-{\theta \over 2}(1+C)<1-{3\theta \over 2}<1
$$

\section{Tools for a Domain Decomposition algorithm  for the Helmholtz equation}
\noindent In the sequel we shall assume {\bf connnectedness} of the open sets $\Omega_m$, $m=1,2$.

\subsection{On the $j$ operator}
We define the linear bounded operator $j$ from $ \Lambda$ to $ \Lambda'$ as the composition of the bounded linear injections
$$ \Lambda \subset H^{1\over 2}_0(\Gamma) \subset L^2(\Gamma)\subset H^{-{1\over 2}}(\Gamma) \subset  \Lambda' 
$$
It will play a key role for Helmholtz equations. Its properties are summarized by:
\bpr  \label{X17}.

1- $j$ is  a compact and one-to-one operator

2- For any $\lambda\in \Lambda$, $j(\overline{\lambda})=\overline{j(\lambda)}$

3- For $m=1,2$  $$\forall \lambda\in \Lambda, \sigma\in \Lambda, (j(\lambda), \sigma)_{\Lambda\Lambda'}=(\sigma,S_m^{-1}j(\overline{\lambda}))_{\Lambda_m}=(\sigma,\overline{\lambda})_{L^2(\Gamma)}$$
$$\forall \lambda\in \Lambda,\quad  (j(\lambda), \overline{\lambda})_{\Lambda,\Lambda'}=\Vert\lambda\Vert_{L^2(\Gamma)}^2
$$

4- $j'=j$

5- For $m=1,2$ the operator $S_m^{-1}j$ is selfadjoint in $\Lambda_m$

6- For $m=1,2$ the operator $jS_m^{-1}$ is selfadjoint in $\Lambda_m'$
\epr

\noindent {\it proof:} Item 1 comes from the Rellich compactness of the injection from $H^{1\over 2}_0(\Gamma)$ to $ L^2(\Gamma)$. 

\noindent Item 3 comes from:

$ (\sigma,S_m^{-1}j(\overline{\lambda}))_{\Lambda_m}=\int_{\Omega_m}\nabla u^\sigma \nabla u^{S_m^{-1}}j(\lambda)=\int_\Gamma \sigma {\partial u^{S_m^{-1}j(\lambda)}\over \partial n_m}=\int_\Gamma \sigma \lambda=(\sigma,\overline{\lambda})_{L^2(\Omega_m)} $

\noindent Item 4 follows item 3 because
$(j(\lambda), \sigma)_{\Lambda \Lambda'}=\int_\Gamma \sigma \lambda=(j(\sigma), \lambda)_{\Lambda\Lambda'}$

\noindent Item 5 follows items 3 and 4 because

$(\sigma,S_m^{-1}j(\overline{\lambda}))_{\Lambda_m}=(j(\lambda), \sigma)_{\Lambda \Lambda'}=(j(\sigma),\lambda )_{\Lambda \Lambda'}=(S_m^{-1}j(\sigma),\overline{\lambda})_{\Lambda_m}$

\noindent Item 6 follows item 5 because

$(jS_m^{-1}\mu,\nu)_{\Lambda_m'}=(S_m^{-1}jS_m^{-1}\mu,S_m^{-1}\nu)_{\Lambda_m}=(S_m^{-1}\mu,S_m^{-1}jS_m^{-1}\nu)_{\Lambda_m}=(\mu,jS_m^{-1}\nu)_{\Lambda_m'}
$

\subsection{The spectrum of the local Helmholtz problems}

\noindent This paragraph is devoted to the study of the operators
$$m=1,2\quad (S_m+ i\gamma j):\Lambda\longrightarrow \Lambda'$$
 related to the Laplace equation, and to the like $(S_m^k+ i\gamma j)$ operators related to the Helmholtz equation. Let $\gamma$ denote a  real number.

\bpr  \label{X18}.

1- $\forall \lambda \in \Lambda\quad ((S_m+ i\gamma j)\lambda,\overline{\lambda})_{\Lambda \Lambda'}=\Vert \lambda\Vert_{\Lambda_m}^2+ i\gamma \Vert \lambda\Vert_{L^2(\Gamma)}$

2- $(S_m+ i\gamma j)$ has a bouded inverse.

3- $\forall \lambda \in \Lambda\quad (S_m+ i\gamma j) \overline{\lambda}= \overline{(S_m- i\gamma j) {\lambda}}\quad and \quad T_m^{ \gamma} \overline{\lambda}= \overline{T_m^{- \gamma}(\lambda)}$
\epr
\noindent {\it proof:}

\noindent 1- is straightforward applying proposition \ref{X17} and notation \ref{X9}

\noindent 2- item 1 shows that $Ker (S_m+ i\gamma j)=\{0\}$ and $Im (S_m+ i\gamma j)$ is closed in $\Lambda'$. It remains to show that $Im (S_m+ i\gamma j)$ is everywhere dense in  $\Lambda'$. By propositions \ref{X10} and \ref{X17}
$$\forall \lambda \in \Lambda, ((S_m+ i\gamma j)\lambda,\eta)_{\Lambda\Lambda'}=0\Longrightarrow
$$
$$\forall \lambda,\,\, ((S_m+ i\gamma j)\eta,\lambda)_{\Lambda\Lambda'}=0\Longrightarrow(S_m+ i\gamma j)\eta=0\Longrightarrow \eta=0
$$
because $(S_m+ i\gamma j)$ is one to one.

\noindent 3- is straightforward.

\bpr  \label{X19} For $m=1,2$ and any $\nu \in \Lambda'$ there exists a unique $u\in H_m$ such that
$$\Delta u=0\,\,{\rm in}\,\,\Omega_m\quad {\rm and} \quad {\partial u\over \partial n_m} + i\gamma j\rho^\Gamma_mu=\nu$$
In fact $u=u^{(S_m+i\gamma j)^{-1}\nu}$ and 
$$\Vert  u\Vert_{H_m}\leq C\Vert \nu\Vert_{\Lambda'}$$
\epr
 \noindent {\it proof:} Uniqueness is straightforward by the Green formula. For existence we apply the previous proposition to get $ \lambda \in \Lambda$ such that $(S_m+ i\gamma j)\lambda=\nu$, and check that  $u=u^\lambda$ solves the problem. 
 
 \medskip
 \noindent The following remark will be crucial to prove convergence of domain decomposition algorithms for the Helmholtz equation:
 
 \brm  \label{X20} with the notation of the preceeding proposition, if $\nu \in L^2(\Gamma)$ then
 $$  {\partial u\over \partial n_m}\in L^2(\Gamma)
 \quad {\rm and}\quad
 \Vert {\partial u\over \partial n_m}\Vert_{L^2(\Gamma)}  \leq C\Vert \nu\Vert_{L^2(\Gamma)}
 $$
 \erm
  \noindent {\it proof:} $ {\partial u\over \partial n_m}=\nu- i\gamma j\rho^\Gamma_mu \in L^2(\Gamma)+\Lambda \subset L^2(\Gamma)$  and  
  
  $$ \Vert {\partial u\over \partial n_m}\Vert_{L^2(\Gamma)}\leq\Vert \nu\Vert_{L^2(\Gamma)}+\vert \gamma\vert \Vert \rho^\Gamma_mu\Vert_{L^2(\Gamma)}\leq  \Vert \nu\Vert_{L^2(\Gamma)}+\vert\gamma \vert \Vert u\Vert_{H_m}$$
  $$\leq  \Vert \nu\Vert_{L^2(\Gamma)}+C\vert \gamma \vert \Vert \nu\Vert_{\Lambda'} \leq   C\Vert \nu\Vert_{L^2(\Gamma)}$$

\bpr \label{X21} For $m=1,2$ and any $f \in L^2(\Omega_m)$ there exists a unique $u\in H_m$ such that
$$\Delta u=f\,\,{\rm in}\,\,\Omega_m\quad {\rm and} \quad {\partial u\over \partial n_m} + i\gamma j\rho^\Gamma_mu=0$$
 and we have the estimate
$$\Vert  u\Vert_{H_m}\leq C\Vert f\Vert_{L^2(\Omega_m)}$$
\epr
  \noindent {\it proof:} let $v\in H^1_0(\Omega_m)$ solve  $\Delta v=f$. Let $w=u-v$. The function $w$ fulfills $\Delta w=0 $ and ${\partial w\over \partial n_m} + i\gamma j\rho^\Gamma_m w={\partial v\over \partial n_m}$. Remark \ref{X29} shows that ${\partial v\over \partial n_m}\in \Lambda'$. The result follows from the previous proposition and the estimate follows  estimates in remark \ref{X29} and proposition  \ref{X19}.

\brm  \label{X23}Using the notations of proposition \ref{X21} we have,  as in remark \ref{X20},
 ${\partial u\over \partial n_m}\in L^2(\Gamma)$ and  the estimate:
$$ \Vert {\partial u\over \partial n_m}\Vert_{L^2(\Gamma)}\leq  C\Vert f\Vert_{L^2(\Omega_m)}$$
\erm
 We can now proceed to compute the eigenfrequencies of the local Helmholtz problems involved in the Domain Decomposition algorithm. For that sake, we will use the following
 \bnt  \label{X24} For $m=1,2$ we denote by $$D_m^\gamma:L^2(\Omega_m)\longrightarrow L^2(\Omega_m)$$
 $$f\longrightarrow u$$ where $u$ is given by proposition \ref{X21}. \ent
 This map has the following properties:
 \bpr  \label{X25}.
 
 1- $D_m^\gamma$ is a compact operator in $L^2(\Omega_m)$
 
 2- the adjoint map of $D_m^\gamma$ for the $L^2(\Omega_m)$ scalar product is $D_m^{-\gamma}$
 
 3- we have $D_m^\gamma \overline{g}= \overline{D_m^{-\gamma} g}$ and $D_m^\gamma D_m^{-\gamma}\overline{f}=\overline{D_m^{-\gamma}D_m^\gamma f}$
 
 4- $Im D_m^\gamma \subset H_m$
 
 5- ${\partial \over \partial n_m}D_m^\gamma f\in L^2(\Gamma)$ 
 \epr
 \noindent {\it proof:}
 
  \noindent 1- translates Rellich compactness of the imbedding of $H_m$ in  $L^2(\Omega_m)$
 
 \noindent  2- If $u=D_m^\gamma f$ and $v=D_m^{-\gamma}g$ then Green formula gives
 $$\int_{\Omega_m}u\overline{g}dx - \int_{\Omega_m}\overline{v}fdx=\int_\Gamma{\partial \overline{v}\over \partial n_m}ud\sigma-\int_\Gamma{\partial u\over \partial n_m}\overline{v}d\sigma=-i\gamma\int_\Gamma u \overline{v}d\sigma+i\gamma\int_\Gamma u \overline{v}d\sigma=0
 $$
 \noindent 3- is straightforward
  
 \noindent  4- follows proposition \ref{X21}
 
  \noindent  5- follows remark \ref{X23}
  \medskip
  
  \noindent We collect the spectral properties of $D^\gamma_m$ in the following
  
  \bpr \label{X26} We denote by $\sigma$ the spectrum of an operator and by $\sigma_p$ the set of its eigenvalues. We have for $m=1,2$:
  \medskip
  
 \noindent 1- $\sigma(D^\gamma_m)=\{0\}\cup \sigma_p(D^\gamma_m)$
  
 \noindent 2- $\mu \in \sigma(D^\gamma_m)\Longleftrightarrow \overline{\mu} \in \sigma(D^{-\gamma}_m)$
 
  \noindent 3- For any $f\in L^2(\Omega_m)$, if $u=D^\gamma_mf$ then
  $$(D^\gamma_m f,f)_{L^2(\Omega_m)}=-\int_{\Omega_m}\vert \nabla u\vert^2dx + i\gamma\int_\Gamma\vert \rho_m^\Gamma u\vert^2d\sigma
  $$
    
 \noindent 4- there exists a constant $c$ such that If $\mu \in \sigma(D^\gamma_m)$, $\mu \neq 0$, then
 $$Re\, \mu\leq 0, \gamma Im\, \mu \in  R_+, {\vert Im \,\mu\vert  \leq c\vert \gamma \vert \,\, \vert Re\, \mu \vert} $$
      
 \noindent 5- If $\gamma \neq 0$ then $$  \sigma(D^\gamma_m)\cap R = \{ 0\}$$
   \epr
  
   \noindent {\it proof:} 
   
   \noindent 1- follows compactness of $D^\gamma_m$ asserted in proposition \ref{X25}.
  
 \noindent 2- is obvious by taking the complex conjugate of the eigenfunction associated with $\mu$.
 
   \noindent 3- $$(D^\gamma_m f,f)_{L^2(\Omega_m)}=\int_{\Omega_m}u\overline{f}dx=\int_{\Omega_m}u\overline{\Delta u}dx=$$
   $$-\int_{\Omega_m}\vert \nabla u\vert^2dx +\int_\Gamma \rho_m^\Gamma u\,\,\overline{{\partial u\over \partial n_m}_{\vert \Gamma}}d\sigma=-\int_{\Omega_m}\vert \nabla u\vert^2dx + i\gamma\int_\Gamma\vert \rho_m^\Gamma u\vert^2d\sigma$$
 
  \noindent 4- if $\mu $ is an eigenvalue of $D^\gamma_m$ with associated eigenfunction $f$, and $u=D^\gamma_m f$, then 
  $$\mu \int_{L^2(\Omega_m)}\vert f\vert^2 dx=-\int_{\Omega_m}\vert u\vert^2dx + i\gamma\int_\Gamma\vert \rho_m^\Gamma u\vert^2d\sigma
  $$ 
  and the result follows with $c$ the constant of continuity of the trace operator $\rho^\Gamma_m$ from $H_m$ to  $L^2(\Gamma)$.
 
   \noindent 5- If $\mu \neq 0$ is a real eigenvalue of $D^\gamma_m$ with associated eigenfunction $f$, then the previous formula shows that $\rho^\Gamma_m u=0$ and, because $u\in Im(D^\gamma_m)$, ${\partial u\over \partial n_m}_{\vert \Gamma}=- i\gamma \rho^\Gamma_m u=0$. On the other hand, the equality $D^\gamma_m f=\mu f$ translates to $\Delta u=  {1\over \mu} u$. Because the Laplace operator is hyperbolic in the direction $n_m$, and both  data on $\Gamma$ are zero, this implies $u=0$ on a neighbourhood of $\Gamma$. Solutions of elliptic equations being analytic, and $\Omega_m$ being connected, this implies $u=0$ on $\Omega_m$. Then $f=0$, which contradicts the assumption on $f$ as an eigenfunction.

 \bth  \label{X28} Let $k^2\in R, k^2\neq 0$.  Let $\gamma \in R, \gamma \neq 0$. For $m=1,2$:  
\medskip
 
 \noindent 1- for any $f$ in $L^2(\Omega_m)$ and $\nu \in \Lambda'$, there exits a unique $u\in H_m$ such that:
 $$\Delta u +k^2u=f\quad  {\partial u\over \partial n_m}- i\gamma j\rho^\Gamma_mu=\nu$$
 
  \noindent 2- we have the estimate
  $$\Vert u \Vert_{H_m}   \leq C(\Vert f\Vert_{L^2(\Omega_m)}+\Vert \nu\Vert_ {\Lambda'})$$

   \noindent 3- if $\nu \in L^2(\Gamma)$ then ${\partial u\over \partial n_m}\in L^2(\Gamma)$ and we have the estimate:

   $$\Vert {\partial u\over \partial n_m} \Vert_{L^2(\Gamma)}  \leq C(\Vert f\Vert_{L^2(\Omega_m)}+\Vert \nu\Vert_ {L^2(\Gamma)})
   $$
   \eth
   \noindent {\it proof:} 
   
   \noindent 1- Let $\lambda=(S_m-i\gamma j)^{-1}\nu$ and let $v=u-u^\lambda$. Then:
   $$\Delta u +k^2u=f\quad {\rm and}\quad  {\partial u\over \partial n_m}- i\gamma j\rho^\Gamma_mu=\nu\Longleftrightarrow$$
   $$ \Delta v =f-k^2u^\lambda -k^2v\quad   {\rm and}\quad {\partial v\over \partial n_m}- i\gamma j\rho^\Gamma_mv=0 \Longleftrightarrow$$
   $$v=D^\gamma_m(f-k^2u^\lambda -k^2v) \Longleftrightarrow
   D^\gamma_m v+k^{-2}v=k^{-2}D^\gamma_m(f-k^2u^\lambda ) $$
   The previous proposition \ref{X26} shows that $-k^{-2}\notin \sigma (D^\gamma_m)$ so this problem is well-posed   
   
   \noindent 2- we have the estimate:
   $$\Vert v\Vert_{L^2(\Omega)}
   \leq C  (k^{-2}\Vert D^\gamma_m f\Vert_{L^2(\Omega)}+\Vert D^\gamma_m u^\lambda \Vert_{L^2(\Omega)})
   \leq C (k^{-2}\Vert  f\Vert_{L^2(\Omega)}+\Vert  u^\lambda \Vert_{L^2(\Omega)})
 $$
 so 
 $$ \Vert u\Vert_{L^2(\Omega)}\leq C' (k^{-2}\Vert  f\Vert_{L^2(\Omega)}+\Vert  u^\lambda \Vert_{L^2(\Omega)})
\leq C' (k^{-2}\Vert  f\Vert_{L^2(\Omega)}+\Vert  u^\lambda \Vert_{H_m})
$$
$$
 \leq C' (k^{-2}\Vert  f\Vert_{L^2(\Omega)}+\Vert  \lambda \Vert_{\Lambda})
\leq C' (k^{-2}\Vert  f\Vert_{L^2(\Omega)}+\Vert  \nu \Vert_{\Lambda'})
 $$
We write
 $$\Delta u =f-k^2u\quad  {\partial u\over \partial n_m}_{\vert \Gamma}- i\gamma j\rho^\Gamma_mu=\nu$$
 and the $H_m$ estimate follows  proposition \ref{X19} and   proposition \ref{X21}.
 
    \noindent 3- If $\nu \in L^2(\gamma)$ we write again 
 $$\Delta u =f-k^2u\quad  {\partial u\over \partial n_m}_{\vert \Gamma}- i\gamma j\rho^\Gamma_mu=\nu$$
 and the estimate follows remark \ref{X20} and remark \ref{X23}.
 \medskip

\subsection{Despres operators and the energy fluxes}

\medskip

We now define the building blocks of the intertwinning operator on the 
fictitious boundary $\Gamma$:
the
Despres operators.
\medskip

\bdf \label{X54}: Let $k\in R$, $\gamma \neq 0$, $\gamma \in R$. For any $\nu \in \Lambda'$, let $u \in H_m$ be the unique solution, given by Theorem \ref{X28}, of the following equation on $\Omega_m$, $m=1,2$:
$$
 \Delta u +k^2 u=0,  \quad
 {\partial u \over \partial n_m }-i\gamma j\rho_m^\Gamma u=\nu \quad {\rm on}\quad \Gamma
$$
 Let $\tilde P^\gamma_m$ be the 
linear bounded 
operator in $\Lambda'$ defined by:
$$
\tilde P^\gamma_m \nu={\partial u \over \partial n_m}\,\,+i\gamma    j\rho_m^\Gamma u\quad {on} \quad \Gamma
$$
\edf

\brm \label{X31} Boundedness of $\tilde P^\gamma_m$ follows Remark \ref{X24} and Proposition \ref{X17}.
\erm
\medskip
\bpr :\label{X30} We obviously have $\tilde P^\gamma_m \tilde P^{-\gamma}_m=\tilde P^{-\gamma}_m\tilde P^\gamma_m=Id_{\Lambda'}$
\epr

\bnt\label{X32} Let  $\tilde A^\gamma$ denote the linear bounded operator in $\Lambda'\times \Lambda'$ given by:
$$ \tilde  A^\gamma= \left( \begin{array}{cc}
0 & -\tilde P^\gamma_1\\
-\tilde P^\gamma_2 &0
\end{array}\right) $$
\ent

\brm \label{X33}:
The inverse of $ \tilde A^\gamma$    in $\Lambda'\times 
\Lambda'$  is the bounded linear
operator given by:
$$
( \tilde  A^\gamma)^{-1}= \left( \begin{array}{cc}
0 & -\tilde P^{-\gamma}_2\\
-\tilde P^{-\gamma}_1 & 0
\end{array} \right)
$$
\erm
\medskip 
 \noindent In order to use conservation of energy, and to gain compactness, we use Theorem \ref{X28} to  introduce:
 \bnt \label{X34} For $\gamma \neq 0$ and $k\in R$, let the bounded operator in $L^2(\Gamma)$ denoted by $P_m^\gamma$ be the restriction of $\tilde P_m^\gamma$ to $L^2(\Gamma)$. Let  the bounded operator in $L^2(\Gamma)\times L^2(\Gamma)$ denoted by $A^\gamma$ be the restriction of $\tilde A^\gamma$ to $L^2(\Gamma)\times L^2(\Gamma)$.
  \ent
 \noindent Conservation of energy fluxes through $\Gamma$ reads:
 
\bpr \label{X35} Let $\gamma \neq 0$ and $m=1,2$. 

\noindent (i) $P_m^\gamma$ is an isometry in $L^2(\Gamma)$: $$\forall \nu \in L^2(\Gamma),\quad 
\Vert  \nu \Vert_{L^2(\Gamma)}=\Vert P_m^\gamma \nu\Vert_{L^2(\Gamma)}$$

\noindent  (ii) $A^\gamma$  is an isometry in $L^2(\Gamma)\times L^2(\Gamma)$:
$$\forall (\nu,\eta) \in L^2(\Gamma)\times L^2(\Gamma),\quad 
\Vert  (\nu, \eta) \Vert_{L^2(\Gamma\times L^2(\Gamma))}=\Vert A^\gamma (\nu, \eta)\Vert_{L^2(\Gamma)\times L^2(\Gamma)}
$$
\epr
\noindent {\it proof:} 
For $\nu \in L^2(\Gamma)$  let $u\in H_m$ solve by Theorem \ref{X28} the following equation:
$$
\Delta u+k^2u=0\quad {\rm in}\quad \Omega_m, \quad {\partial u\over \partial n_m}-i\gamma j \rho_m^\gamma u=\nu \quad \Gamma
$$
\noindent
 Multiplying by $\overline{ u}$ the equation fulfilled by 
$u$ and 
integrating on $\Omega_m$ gives
$$\int_ {\Omega_m} \vert \nabla u\vert^2dx -k^2\int_ {\Omega_m} \vert  u\vert^2dx=
\int_\Gamma {{\partial u}\over {\partial n_m}}{\overline u}d\sigma$$
\noindent Taking the imaginary part gives $${\cal I}m \int_\Gamma {{\partial u}\over {\partial n_m}}{\overline u}d\sigma =0$$

\noindent The result follows integration of the following
 identity on $\Gamma$:
$$\vert {{\partial u}\over {\partial n_m}}+i\gamma j\rho^\gamma_m u\vert^2-\vert {{\partial u}\over
 {\partial n_m}} -i\gamma j\rho^\gamma_m u\vert^2={ 4i\gamma}\,\,{\cal I}m\,\, \overline{ u}
{{\partial u}\over {\partial n_m}} $$

\noindent An important consequence of this property will be crucial in the next section:
\bcr
\label{X36}: For $m=1,2$ and $\gamma \neq 0$,  

\noindent (i)   $P_m^\gamma$ is a  normal operator in $L^2(\Gamma)$

\noindent (ii)  $A^\gamma$  is a normal operator in $L^2(\Gamma)\times L^2(\Gamma)$
\ecr

\subsection {Spectral properties of the Despres operators}
\medskip
In the preceeding section we proved that the Despres operator  $P_m^\gamma$, ( $m=1,2$ and $\gamma \neq 0$)
is a bijective isometry in $L^2(\Gamma)$, and consequently a normal operator in $L^2(\Gamma)$. It follows that its spectrum 
is a  subset of the unit circle in the complex plane. We now investigate this 
spectrum more accurately.
\bdf \label{X37} Let $\gamma \neq 0$ and $m=1,2$. Let $C_m^\gamma$ be the operator in $L^2(\Gamma)$ given by:
$$\forall \nu \in L^2(\Gamma), \quad C_m^\gamma \nu= j\,\,\rho_m^\Gamma u$$
where $u \in \Lambda$ is the solution given by Theorem \ref{X28} of the equation:
$$(\Delta+k^2)u=0\quad {\rm in}\quad \Omega_m,\quad  {{\partial u}\over {\partial n_m}}-i\gamma \rho^\Gamma_mu=\nu \quad {\rm on} \quad \Gamma$$
\edf

\noindent Compactness of the injection $j$ from $\Lambda$ to $L^2(\gamma)$ gives:

\bpr \label{X38}: 

\noindent (i) For $m=1,2$ and $\gamma \neq 0$, $ P_m^\gamma=I+2i\gamma C_m^\gamma$

\noindent (ii) $C_m^\gamma$ is a normal and compact operator in $L^2(\gamma)$
\epr

\bnt \label{X39} Let 

\noindent  (i) $\,\,\Sigma^\gamma_m$ denote the spectrum of $P^\gamma_m$ 

\noindent  (ii) $\Sigma^{Dir}_m$ denote the sequence of eigenvalues of the Laplace operator on $\Omega_m$ with Dirichlet boundary condition on $\partial \Omega_m$, i.e. 

$-k^2 \notin \Sigma^{Dir}_m$ {if and only if the following problem is well posed:}
$$(\Delta +k^2)u=f \in L^2(\Omega_m), \quad u\in H^1_0(\Omega_m)$$
\ent

\noindent We have:

\bpr \label{X40}: Let $\gamma \neq 0$, and $m=1,2$. 

\noindent (i) $1$ belongs to $\Sigma^\gamma_m$

\noindent  (ii)  $1$ is an eigenvalue of $P^\gamma_m$ if and only if $-k^2 \in \Sigma^{Dir}_m$
\epr
\noindent {\it proof:}

\noindent (i) Because $C^\gamma_m$ is a compact normal operator (Proposition \ref{X38}) in $L^2(\Omega_m)$, its spectrum is a sequence of eigenvalues and its limit zero. This implies that $1$ is in the closure of $\Sigma^\gamma_m$, which is closed.

\noindent (ii) If $-k^2 \in  \Sigma^{Dir}_m$,  then there exist an eigenfunction $\varphi_k$ of the Laplace operator  such that
$$\Delta \varphi_k+k^2\varphi_k=0 \quad {\rm and}\quad { \varphi_k}_{\vert \Gamma}=0$$
Let $\nu_k={\partial  \varphi_k \over \partial n_m}_{\vert \Gamma}$. We have $\nu_k \neq 0$ or else $ { \varphi_k}$ solving an elliptic equation, condition ${ \varphi_k}_{\vert \Gamma}=0$ and connectedness of $\Omega_m$  would imply ${ \varphi_k}=0$ on $\Omega_m$, which contradicts the fact that ${ \varphi_k}$ is an eigenfunction. Obviously $\nu_k$ is an eigenvector of $P^\gamma_m$ for the eigenvalue $1$.

\noindent Conversely,  if $1$ is an eigenvalue of $P^\gamma_m$, with eigenvector $\nu \neq 0$, then there exist $u$ in $H_m$ such that 
$$\Delta u+k^2 u=0 \quad {\rm and}\quad \nu ={\partial u \over \partial n_m}_{\vert \Gamma}-i\gamma j \rho^\Gamma_m u= {\partial u \over \partial n_m}_{\vert \Gamma}+i\gamma j \rho^\Gamma_m u$$
This implies $\rho^\Gamma_m u=0$, which added to the fact that $u \in H_m$ implies $u\in H^1_0(\Omega_m)$. So $u$ is an eigenfunction of the Laplace operator for the eigenvalue $-k^2$, provided $u$ is not identically zero. And this is ruled out because ${\partial u \over \partial n_m}_{\vert \Gamma}=\nu \neq  0$

\noindent The following Theorem gives a complete spectral description of the Despres operators:

\bth \label{X42}: Assume that $-k^2 \notin  \Sigma^{Dir}_m$. For $\gamma \neq 0$ and $m=1,2$:

\noindent (i) $\Sigma^\gamma_m=\{1\}\cup (e^{i\sigma^n_m})_{n\in N}$ where $(\sigma^n_m)_{n\in N}$ is a sequence of numbers with $\sigma^n_m \in R$, $\sigma^n_m \neq 0$, and $\sigma^n_m \longrightarrow 0$ when $n\longrightarrow \infty$.

\noindent (ii) For each $n\in N$, $e^{i\sigma^n_m}$ is an  eigenvalue of $P^\gamma_m$, 
 with finite multiplicity.

\noindent (iii)  $L^2(\Gamma)$ is the Hilbert direct sum of the eigenspaces associated with the eigenvalues $e^{i\sigma^n_m}$ 
\eth
{\it Proof:} $C_m^\gamma$ is a  normal and compact  operator in $L^2(\Gamma)$ (proposition  
\ref{X38}). Its kernel is trivial (proposition  
\ref{X40}). The diagonalization theorem gives a sequence 
$(\lambda_m^n)_{n\in N}$, 
($\lambda_m^n\neq 0$),
 for its eigenvalues (they have finite multiplicity), with limit zero, and $L^2(\Gamma)$  is the Hilbert direct sum of the associated eigenspaces. The diagonalization of $P^\gamma_m=I+2i\gamma C_m^\gamma$ follows, 
with eigenvalues  $(1+2i\gamma \lambda_m^n)_{n\in N}$, $(1+2i\gamma \lambda_m^n\neq 1)$. 
Proposition \ref{X35} implies that these eigenvalues have modulus one: we set 
$1+2i\gamma \lambda_m^n=e^{i\sigma^n_m}$.

\subsection{Spectral Properties of  $P^\gamma_1P^\gamma_2$ and $P^\gamma_2P^\gamma_1$}
\medskip

Properties of the intertwinning operator $A^\gamma$  rely  heavily on the spectral 
properties
of  $P^\gamma_1P^\gamma_2$ and $P^\gamma_2P^\gamma_1$ that we investigate now.
\medskip 

\noindent We first list  obvious properties which follow from the  previous 
section: 

\bpr \label{X41}:
 
\noindent (i) $P^\gamma_1P^\gamma_2$ and $P^\gamma_2P^\gamma_1$ are isometric bijections in $L^2(\gamma)$.

\noindent (ii) $P^\gamma_1P^\gamma_2$ and $P^\gamma_2P^\gamma_1$ are normal operators in $L^2(\gamma)$.

\noindent (iii) $P^\gamma_1P^\gamma_2-I$ and $P^\gamma_2P^\gamma_1-I$ are compact operators in $L^2(\gamma)$.
\epr 

\noindent {\it Proof:}

\noindent (i) follows the fact that $P^\gamma_1$ and $P^\gamma_2$ are   isometric bijections in $L^2(\Gamma)$.

\noindent (ii) follows  (i)

\noindent (iii) follows proposition \ref{X38} through:
$$P^\gamma_1P^\gamma_2=(I+2i\gamma C^\gamma_1)(I+2i\gamma C^\gamma_2)=I+2i\gamma C^\gamma_1+2i\gamma C^\gamma_2-4\gamma^2 C^\gamma_1 C^\gamma_2
$$
\medskip

\noindent An important property that we sall need is the spectral status of 1: 

\bpr \label{X43}: 
 1 is not an eigenvalue of $P^\gamma_1P^\gamma_2$ or $P^\gamma_2P^\gamma_1$ in $L^2(\Gamma)$.
\epr

\noindent {\it Proof:} By symmetry, it is enough to prove it for $P^\gamma_1P^\gamma_2$. Let $\nu \in L^2(\Gamma)$ be such that 
$P^\gamma_1P^\gamma_2\nu=\nu$. This translates to the existence of $u_1 \in H_1$ and $u_2\in H_2$ satisfying:
$$\Delta u_2+k^2u_2=0\,\, {\rm in}\,\, \Omega_2;
\quad {{\partial u_2}\over {\partial n_2}}-iku_2=\nu,\,\, {{\partial u_2}\over 
{\partial n_2}}+iku_2=P^\gamma_2\nu \quad {\rm on} \quad \Gamma $$ 
$$\Delta u_1+k^2u_1=0\,\, {\rm in}\,\, \Omega_1; 
\quad {{\partial u_1}\over {\partial n_1}}-iku_1=P^\gamma_2\nu,\,\, {{\partial u_1}\over
 {\partial n_1}}+iku_1=\nu \quad {\rm on} \quad \Gamma $$
\noindent This implies that on $\Gamma$ these functions fulfill:
$${{\partial u_1}\over {\partial n_1}}-iku_1-{{\partial u_2}\over 
{\partial n_2}}-iku_2=0;\quad{{\partial u_1}\over {\partial n_1}}+iku_1-
{{\partial u_2}\over 
{\partial n_2}}+iku_2=0
$$
\noindent Adding and substracting gives:
$$u_1=-u_2\quad {\rm and}\quad {{\partial u_1}\over {\partial n_1}}=
{{\partial u_2}\over 
{\partial n_2}}
$$
\noindent We define $u$ on $\Omega$ as $u_{\vert \Omega_1}=u_1$ and 
$u_{\vert \Omega_2}=-u_2$. It solves the Helmholtz equation on $\Omega_1$ and 
$\Omega_2$, its has no jump accross $\Gamma$, neither has its normal derivative. 
So it solves Helmholtz equation on $\Omega$. Moreover $u_{\vert \partial 
\Omega}=0$. Assumption {\bf (A) } gives $u=0$, so $u_1=0$ and $u_2=0$; 
and $\nu =0$ follows.
\medskip

\bpr \label{X44}:
 
\noindent  (i) The spectrum of $P^\gamma_1P^\gamma_2$  in $L^2(\Gamma)$ is
 $\{ 1 \}\cup(e^{i\tau_{12}^n})_{n\in N}$,
where $(\tau_{12}^n)_{n\in N}$ is an infinite sequence of real numbers,
 $\tau_{12}^n \neq 0$, and 
$\tau_{12}^n \longrightarrow 0$ when $n\longrightarrow\infty$. 

\noindent  (ii) $(e^{i\tau_{12}^n})_{n\in N}$ is the set of  eigenvalues of $P^\gamma_1P^\gamma_2$.
They have finite multiplicity. If we denote by $E_{12}^n$
the eigenspace associated with $e^{i\tau_{12}^n}$, then $L^2(\Gamma)$ is the Hilbert direct sum of the subspaces
$(E_{12}^n)_{n\in N}$

\noindent  (iii) $P^\gamma_2P^\gamma_1$ has the same properties, and we set the obvious
 notations:  
$(e^{i\tau_{21}^n})_{n\in N}$ for eigenvalues and $(E_{21}^k)_{n\in N}$ for 
eigenspaces.
\epr
\medskip

\noindent  {\it Proof:} The operator $P^\gamma_1P^\gamma_2-I$ is a normal compact operator 
(proposition \ref{X41}). So by the diagonalization theorem its spectrum is the union of $\{0\}$ and 
an infinite sequence of eigenvalues with finite multiplicity
$(t^n)_{n\in N}$, ($t^n\neq 0$). Zero is not an eigenvalue of $P^\gamma_1P^\gamma_2-I$ (proposition 
\ref{X43}). 
So the whole set of eigenvalues is $(t_n)_{n\in N}$. If $E_{12}^n$ 
denotes the eigenspace associated with $t^n$, then $L^2(\Gamma)$ is the Hilbert direct sum of $(E_{12}^n)^{n\in N}$. By proposition \ref{X41} we know that
 $P^\gamma_1P^\gamma_2$ is 
an  isometry in $L^2(\Gamma)$, so $\vert1+ t^n \vert=1$, and we write it: $ t^n 
=e^{i\tau_{12}^n}$. The theorem translates proven properties of $t_n$ 
into properties of $\tau_{12}^n$.
\medskip

\noindent In order to study the relationship between 
$(\tau_{12}^n)_{n\in N}$ and $(\tau_{21}^n)_{n\in N}$, and between 
$(E_{12}^n)_{n\in N}$  and $(E_{21}^n)_{n\in N}$, we prove the 
following lemmi:
\medskip

\blm \label{X45} Let $\gamma \neq 0$. For $m=1,2$
and any $ \nu \in L^2(\Gamma)$: $ \overline{P^\gamma_m \nu}=P^{-\gamma}_m \overline{\nu} $
\elm 

\noindent {\it Proof:} Let $m=1$ or $2$. By definition of $P^\gamma_m$ 
there exists ${u^m}
\in H_m$
 with:
$$\Delta  {u^m}+k^2 {u^m}=0;\quad  {{\partial 
 {u^m}}\over
 {\partial n_m}}_{\vert \Gamma}-i\gamma j \rho^\Gamma_m {u^m}= {\nu}; \quad {{\partial  {u^m}}\over 
{\partial n_m}}_{\vert \Gamma}+i\gamma j \rho^\Gamma_m {u^m}= {P^\gamma_m\nu}
$$
\noindent Taking the complex conjugate of these equalities gives  $ 
\overline{u^m}\in H_m$ such that
$$\Delta  \overline{u^m}+k^2 \overline{u^m}=0;\quad  {{\partial  \overline{u^m}}\over
 {\partial n_m}}_{\vert \Gamma}+i\gamma j \rho^\Gamma_m \overline{u^m} =\overline{\nu}; \quad
 {{\partial  \overline{u^m}}\over 
{\partial n_m}}_{\vert \Gamma}-i\gamma j \rho^\Gamma_m \overline{u^m}= \overline{P^\gamma_m\nu}
$$
\noindent which by definition of $P^\gamma_m$ writes
$$ P^{-\gamma}_m\overline{\nu}=\overline{P^\gamma_m\nu}
$$

\blm \label{X46} Let $\gamma \neq 0$

\noindent (i) If $\lambda$ is an eigenvalue of $P^{\gamma}_1P^{\gamma}_2$ (resp $P^{\gamma}_2P^{\gamma}_1$)
for the eigenvector
 $\nu$ then it  is an eigenvalue  of $P^{\gamma}_2P^{\gamma}_1$ (resp $P^{\gamma}_1P^{\gamma}_2$) with 
associated eigenvector $\overline{\nu}$

\noindent (ii) For all $n\in N$, $\tau_{12}^n=\tau_{21}^n \quad (mod\,\, 2\pi)$; 
we denote it by $\tau_n$

\noindent (iii) If we denote by $\overline{C}$ the set of complex conjugates 
 of distributions in a set $C$, then, for any ${n\in N}$, 
$$\overline{E_{21}^n}=E_{12}^n\quad {\rm and}\quad \overline
{E_{12}^n}=E_{21}^n$$
 \elm

\noindent {\it Proof:} 

\noindent (i) If $P^{\gamma}_1P^{\gamma}_2\nu=\lambda \nu$ then 
$\overline{P^{\gamma}_1P^{\gamma}_2\nu}=\overline{\lambda}\overline{\nu}$ 
which by lemma \ref{X45} writes $P^{-\gamma}_1P^{-\gamma}_2\overline{\nu}=
\overline{\lambda}\overline{\nu}$, which implies, by proposition \ref{X30},  
$P^{\gamma}_2P^{\gamma}_1\overline{\nu}={1\over \overline{\lambda}} \overline{\nu}
={\lambda} \overline{\nu}$ because 
 ${1\over \overline{\lambda}}=\lambda$ by proposition \ref{X44}.  
 
 \noindent (ii) follows (i) and a  renumbering. 

\noindent (iii) follows (i) because it gives:
$$\overline{E_{21}^n}\subset E_{12}^n\quad {\rm and}\quad \overline
{E_{12}^n}\subset E_{21}^n$$
\noindent but then
$$E_{21}^n\subset \overline
{E_{12}^n}\subset E_{21}^k \quad {\rm and}\quad E_{12}^n\subset \overline
{E_{21}^n}\subset E_{12}^n$$
\noindent which gives $$\overline{E_{21}^n}=E_{12}^n\quad {\rm and}\quad
\overline{E_{12}^n}=E_{21}^n$$

\blm \label{X47}: Let $\gamma \neq 0$. For any $n\in N$, 
  $P^\gamma_1E_{21}^n=E_{12}^n $ and $P^\gamma_2E_{12}^n=E_{21}^n $
\elm
\medskip

\noindent {\it Proof:} 
 Let $\nu \neq 0 \in E_{12}^n$ then $P^\gamma_1P^\gamma_2\nu=e^{i\tau_n} \nu$ so 
$P^\gamma_2P^\gamma_1P^\gamma_2\nu=e^{i\tau_n}P^\gamma_2\nu$ which proves that $P^\gamma_2\nu \in E_{21}^n$ 
($P^\gamma_2\nu \neq 0$ because $P^\gamma_2$ is bijective (proposition \ref{X30})). This writes
$P^\gamma_2E_{12}^n\subset E_{21}^n$. Proposition \ref{X44}, lemma \ref{X46} and invertibility of $P^\gamma_2$ 
(proposition \ref{X30}) give
 $dim E_{21}^n=
dim E_{12}^n=dim P^\gamma_2E_{12}^n$ so $P^\gamma_2E_{12}^n=E_{21}^n $.
\medskip

\noindent The following algebraic property and its consequences on the eigenprojectors 
(next theorem) are a key for understanding the geometric properties of 
the intertwinning operator:

\blm \label{X48}: Let $\gamma \neq 0$. For any $\mu \notin \{1\}\cup (e^{i\tau_n})_{n\in N}$:
$$(P^\gamma_1P^\gamma_2-\mu I)^{-1}P^\gamma_1=P^\gamma_1(P^\gamma_2P^\gamma_1-\mu I)^{-1}
$$
$$
(P^\gamma_2P^\gamma_1-\mu I)^{-1}P^\gamma_2=P^\gamma_2(P^\gamma_1P^\gamma_2-\mu I)^{-1}
$$
\elm
\medskip

\noindent {\it Proof:} we prove the second assertion, using resolvant identity:
$$ (P^\gamma_2P^\gamma_1-\mu I)^{-1}P^\gamma_2=(P^\gamma_2P^\gamma_1-\mu I)^{-1}P^\gamma_2(P^\gamma_1P^\gamma_2-\mu I)
(P^\gamma_1P^\gamma_2-\mu I)^{-1}=$$
$$[[I+\mu(P^\gamma_2P^\gamma_1-\mu I)^{-1}]P^\gamma_2-\mu (P^\gamma_2P^\gamma_1-\mu I)^{-1}P^\gamma_2](P^\gamma_1P^\gamma_2-\mu I)^{-1}
=P^\gamma_2(P^\gamma_1P^\gamma_2-\mu I)^{-1}
$$

\blm \label{X49} Let $\gamma \neq 0$. For any $n\in N$, if $\Pi_{12}^n$ (resp. $\Pi_{21}^n$)denotes 
the spectral projector of the operator $P^\gamma_1P^\gamma_2$ (resp. $P^\gamma_2P^\gamma_1$) 
on the eigenspace $E_{12}^n$ (resp. $E_{21}^n$) then we have
$$P^\gamma_1\Pi_{21}^n=\Pi_{12}^nP^\gamma_1 \quad {\rm and}\quad P^\gamma_2\Pi_{12}^n=\Pi_{21}^nP^\gamma_2$$
\elm
\medskip

\noindent {\it Proof:} By symmetry it is enough to prove the first formula. 
Let $C_n$ denote a positively oriented curve in the 
complex plane, which winds one time around the eigenvalue $e^{i\tau_n}$, 
and none around any other eigenvalue, then the Dunford integral representation 
formula gives:
$$\Pi_{12}^n ={-1\over 2i\pi}\int_{C_n} (P^\gamma_1P^\gamma_2-\mu I)^{-1}d\mu \quad {\rm and}  
\quad  \Pi_{21}^n ={-1\over 2i\pi}\int_{C_n} (P^\gamma_2P^\gamma_1-\mu I)^{-1}d\mu    $$
\noindent The previous lemma gives the following:
$$P^\gamma_2\Pi_{12}^n = {-1\over 2i\pi}\int_{C_n} P^\gamma_2(P^\gamma_1P^\gamma_2-\mu I)^{-1}d\mu = 
{-1\over 2i\pi}\int_{C_n} (P^\gamma_2P^\gamma_1-\mu I)^{-1}P^\gamma_2d\mu=\Pi_{21}^nP^\gamma_2 $$

\subsection{Spectral Properties of $A^\gamma$}
We recall that
$$A^\gamma=\left( \begin{array}{cc}
0 & -P^\gamma_1\\
-P^\gamma_2 & 0 
\end{array}\right)
$$
\noindent and that this operator is a bijective isometry 
in $L^2(\Gamma)\times L^2(\Gamma)$\
\medskip

\bth \label{X50}:  Let $\gamma \neq 0$. 

\noindent (i) If $\lambda \notin \{\pm1\}\cup
 (\pm e^{i{\tau_n \over 2}})_{n \in N}$ then $\lambda$ belongs to the resolvant 
set of $A^\gamma$ and
$$(A^\gamma-\lambda I)^{-1}=\left( \begin{array}{cc}
 {\lambda}(P^\gamma_1P^\gamma_2-\lambda^2 I)^{-1}& -(P^\gamma_1P^\gamma_2-\lambda^2 I)^{-1}P^\gamma_1\\
  -(P^\gamma_2P^\gamma_1-\lambda^2 I)^{-1}P^\gamma_2 & {\lambda}(P^\gamma_2P^\gamma_1-\lambda^2 I)^{-1}
\end{array}\right)$$

\noindent (ii) For any $n\in N$, $\pm e^{i{{\tau_n} \over 2}}$ is an eigenvalue 
of $A^\gamma$ with  associated eigenspace: 
\beq  F_n^\pm=\{ (\mu \,\, , \,\,\mp e^{-{{i\tau_n} \over 2}}P^\gamma_2\mu);
\quad \mu \in E^n_{12} \}
\eeq
and  associated eigenprojector :
\beq  P_n^\pm = \left( \begin{array}{cc}
{1\over 2}\Pi_{12}^n &  \mp {1\over 2}e^{-i{{\tau_n}\over {2}}} P^\gamma_1 \Pi_{21}^n\\
 \mp {1\over 2}e^{-i{{\tau_n}\over {2}}}P^\gamma_2 \Pi_{12}^n  & {1\over 2}\Pi_{21}^n
\end{array}
\right)
\eeq
\noindent (iii) $ \{\pm 1\}$  belong to the spectrum of $A^\gamma$ and  are  not  eigenvalues of $A^\gamma$

\noindent (iv) $(F_n^\pm)_{n\in N;\pm}$ is an orthogonal family of subspaces and we have the Hilbert decomposition
$$L^2(\Gamma)\times L^2(\Gamma)= ( \oplus_0^\infty  F_n^+) \oplus  (\oplus_0^\infty  F_n^-)$$

\noindent (v) The following series are strongly convergent in ${\cal L}(L^2(\Gamma)\times L^2(\Gamma))$:
$$I=\sum_0^\infty P_n^+ + \sum_0^\infty P_n^-$$
$$A^\gamma=\sum_0^\infty e^{i{\tau_n \over 2}}P_n^+ -\sum_0^\infty e^{i{\tau_n \over 2}}P_n^-$$
\eth
\medskip

\noindent {\it Proof:} 

\noindent (i) let $\lambda \notin \{\pm1\}\cup
 (\pm e^{i{\tau_n \over 2}})_{n \in N}$

\noindent {\it $(A^\gamma-\lambda I)$ is injective:} let $(\varphi , \psi)\in 
L^2(\Gamma)\times L^2(\Gamma)$ be such that 
$$(A^\gamma-\lambda I)\left( \begin{array}{c}\varphi \\ \psi \end{array} \right)
=\left( \begin{array}{c}
0 \\ 0 \end{array}\right)$$
\noindent this writes
$$P^\gamma_1\psi +\lambda \varphi=0 \quad {\rm and}\quad P^\gamma_2\varphi +\lambda \psi =0$$
\noindent which implies
$$P^\gamma_2P^\gamma_1\psi -\lambda^2 \psi=0 \quad {\rm and}\quad P^\gamma_1P^\gamma_2\varphi -\lambda^2 
\varphi =0$$
\noindent This implies $\varphi=\psi=0$ by proposition \ref{X44} and lemma \ref{X46}.
\medskip

\noindent {\it $(A^\gamma-\lambda I)$ is surjective:}  For any $(\xi , \eta)\in 
L^2(\Gamma) \times L^2(\Gamma)$ let:
$$\varphi = (P^\gamma_1P^\gamma_2-\lambda^2 I)^{-1}({\lambda}\xi -P^\gamma_1\eta)\quad {\rm and}
 \quad \psi = (P^\gamma_2P^\gamma_1-\lambda^2 I)^{-1}({\lambda}\eta -P^\gamma_2\xi)
$$

\noindent We have, by lemma \ref{X48}
$$P^\gamma_2\varphi+\lambda \psi=P^\gamma_2(P^\gamma_1P^\gamma_2-\lambda^2 I)^{-1}({\lambda}\xi -P^\gamma_1\eta)
+\lambda (P^\gamma_2P^\gamma_1-\lambda^2 I)^{-1}({\lambda}\eta -P^\gamma_2\xi)
$$
$$=(P^\gamma_2P^\gamma_1-\lambda^2 I)^{-1}({\lambda}P^\gamma_2\xi -P^\gamma_2P^\gamma_1\eta)
+\lambda (P^\gamma_2P^\gamma_1-\lambda^2 I)^{-1}({\lambda}\eta -P^\gamma_2\xi)
$$
$$=(P^\gamma_2P^\gamma_1-\lambda^2 I)^{-1}(-P^\gamma_2P^\gamma_1\eta{+\lambda^2}\eta)=-\eta
$$

\noindent Similarly
$$P^\gamma_1\psi+\lambda \varphi=-\xi
$$
\noindent These  two equalities write
$$(A^\gamma-\lambda I)\left( \begin {array}{c}\varphi \\ \psi \end{array}\right)
=\left( \begin{array}{c}\xi \\ \eta \end{array} \right)
$$
\noindent So surjectivity is proven. These expressions for $(\varphi, \psi)$ give the formula for the resolvent of $A$.
\medskip

\noindent (ii) By definition of $E^n_{12}$ we have for any $\mu \in E^n_{12}$:
$$A^\gamma(\mu \,\, , \,\,\mp e^{-i{{\tau_n} \over 2}}P^\gamma_2\mu)=
(\pm e^{-i{{\tau_n} \over 2}}P^\gamma_1P^\gamma_2\mu \,\, , \,\,-P^\gamma_2\mu )=$$
$$(\pm e^{i{{\tau_n} \over 2}}\mu \,\, , \,\, -P^\gamma_2\mu )=
\pm e^{i{{\tau_n} \over 2}}
(\mu \,\, , \,\, \mp e^{-i{{\tau_n}\over 2}}P^\gamma_2\mu )$$

\noindent Because $E^n_{12}\neq \{ 0\}$, this proves  that $\pm e^{i{{\tau_k} \over 2}} $ is an eigenvalue of $A^\gamma$.
This proves moreover that $F_n^\pm$ is a subset of the eigenspace of $A^\gamma$ 
associated with the eigenvalue $\pm e^{i{{\tau_n} \over 2}} $.

\noindent On the other hand, if $(\xi ,\eta)$ is an eigenvector of $A^\gamma$ for 
the eigenvalue $\pm e^{i{{\tau_n} \over 2}} $ then 
$$-P^\gamma_1\eta=\pm e^{i{{\tau_n} \over 2}}\xi \quad {\rm and}\quad -P^\gamma_2\xi=
\pm e^{i{{\tau_n} \over 2}}\eta
$$
\noindent This implies 
$$P^\gamma_1P^\gamma_2\xi=\mp e^{i{{\tau_n} \over 2}}P^\gamma_1\eta=e^{i{\tau_n} }\xi
\quad {\rm so}\quad \xi \in E^n_{12}
$$
\noindent and
$$\eta =\mp e^{-i{{\tau_n} \over 2}}P^\gamma_2\xi
$$
 \noindent This completes the caracterisation of the eigenspace. 
 \medskip

 \noindent We compute now the eigenprojector: 
for this sake, we make a choice of a branch for $\sqrt {z}$. 
We take  a positively oriented curve  $C_n^\pm$  in the 
complex plane which winds one time around $ \pm e^{i{{\tau_n}\over {2}}}$ 
and not around  
$ \mp e^{i{{\tau_n}\over {2}}}$ nor
does it wind around  any 
$\pm e^{i{{\tau_{n'}}\over {2}}}$ for $n'\neq n$. Let $D_n$ be the image of 
$C_n^\pm$ by the function $z\rightarrow z^2$. $D_n$ winds one time around 
$ e^{i{{\tau_n}}}$ and does not wind around $ e^{i{{\tau_{n'}}}}$ for 
$n'\neq n$. Let $D_n'$ wind one time around 
$ e^{i{{\tau_n}}}$, lying in the interior set delimited by $D_n$.
\medskip

\noindent The eigenprojector is given by the Dunford formula:
$$P_n^\pm = {-1\over 2i\pi }\int_{C_n^\pm}(A^\gamma-\lambda I)^{-1}d\lambda$$
Using the representation formula given by (i) for $(A^\gamma-\lambda I)^{-1}$ leads to compute 
 integrals of two different types:

 \noindent For the first type it is straightforward and gives:
$${-1\over 2i\pi }\int_{C_n^\pm}{\lambda}(P^\gamma_1P^\gamma_2-\lambda^2 I)^{-1}d\lambda 
= {-1\over 4i\pi }\int_{D_n}(P^\gamma_1P^\gamma_2-\lambda I)^{-1}d\lambda 
= {1\over 2}\Pi_{12}^n$$

 \noindent For the second type, we first use the resolvant identity to have:
$${-1\over 2i\pi }\int_{D_n}(P^\gamma_1P^\gamma_2-\lambda I)^{-1}\Pi^n_{12}{d\lambda \over 2\sqrt{\lambda}}=
$$
$${-1\over 2i\pi }\int_{D_n}(P^\gamma_1P^\gamma_2-\lambda I)^{-1}{d\lambda \over 2\sqrt{\lambda}}\,\,
{-1\over 2i\pi }\int_{D_n'}(P^\gamma_1P^\gamma_2-\mu I)^{-1}d\mu =
$$
$${-1\over 2i\pi }{-1\over 2i\pi }\int_{D_n}\int_{D_n'}(P^\gamma_1P^\gamma_2-\lambda I)^{-1}(P^\gamma_1P^\gamma_2-\mu I)^{-1}
{{d\lambda d\mu}\over 2\sqrt{\lambda}}=
$$
$${-1\over 2i\pi }{-1\over 2i\pi }\int_{D_n}(P^\gamma_1P^\gamma_2-\lambda I)^{-1}\left( \int_{D_n'}{{d\mu }\over {\lambda -\mu}} \right) {{d\lambda }\over 2\sqrt{\lambda}}
$$
$$
-{-1\over 2i\pi }{-1\over 2i\pi }\int_{D_n'}(P^\gamma_1P^\gamma_2-\mu I)^{-1}\left( \int_{D_n}{{d\lambda }\over {2\sqrt{\lambda} ({\lambda -\mu}}) }\right)d\mu=
$$
$${-1\over 2i\pi }\int_{D_n'}(P^\gamma_1P^\gamma_2-\mu I)^{-1}{{d\mu }\over {2\sqrt{\mu} } }
$$

\noindent We compute now the second type of integral using this equality,   properties of $\Pi^n_{12}$, and  lemma \ref{X49} to have:
$$ -{-1\over 2i\pi }\int_{C_n^\pm}(P^\gamma_1P^\gamma_2-\lambda^2 I)^{-1}P^\gamma_1d\lambda
=-\left( {-1\over 2i\pi }\int_{D_n}(P^\gamma_1P^\gamma_2-\lambda I)^{-1}{{d\lambda}\over 
{2\sqrt{\lambda}}}\right) P^\gamma_1=
$$
$$-\left( {-1\over 2i\pi }\int_{D_n}(P^\gamma_1P^\gamma_2-\lambda I)^{-1}\Pi^n_{12}{{d\lambda}\over 
{2\sqrt{\lambda}}}\right) P^\gamma_1=-\left( {-1\over 2i\pi }\int_{D_n}{1\over {(e^{i\tau_n}-\lambda I)}}\Pi^n_{12}{{d\lambda}\over 
{2\sqrt{\lambda}}}\right) P^\gamma_1
$$
$$
= \mp {1\over 2}e^{-i{{\tau_n}\over {2}}}  \Pi_{12}^n  P^\gamma_1=\mp {1\over 2}e^{-i{{\tau_n}\over {2}}} P^\gamma_1 \Pi_{21}^n  
$$
\medskip

\noindent (iii) $\pm 1$ are limits of the sequence of eigenvalues $(\pm e^{-i{\tau_n}\over 2})_{n\in N}$ so they belong to the spectrum of $A$. These values are not eigenvalues, or else $1$ is an eigenvalue of $
P^\gamma_2 P^\gamma_1$ and $P^\gamma_1 P^\gamma_2$, which is ruled out by proposition \ref{X43}
\medskip

\noindent (iv) and (v) Assertions (i), (ii) and (iii) prove that the spectrum of $A^\gamma$ is $ \{\pm 1\}\cup  (\pm e^{i{{\tau_n}\over {2}}})_{n\in N}$. Normality of $A^\gamma$ (corollary \ref{X36})
implies orthogonality of the family $(F_n^\pm)_{n\in N}$, and gives the decomposition of $I$ and $A^\gamma$ as series of the eigenprojectors $(P^\pm_n)$

\brm : Notice that  the expression of $F_n^\pm$ in (ii) of the previous proposition  is symmetric:
in fact we have
$$\{ (\mu \,\, ,\,\, \mp e^{-i{{\tau_n}\over 2}}P^\gamma_2\mu);\quad \mu \in E^n_{12}\}
=
\{ (\mp e^{-i{{\tau_n}\over 2}}P^\gamma_1\mu'\,\, ,\,\,\mu' );\quad \mu' \in E^n_{21}\}
$$
 \erm

\noindent  this is because $P^\gamma_1E_{21}^n=E_{12}^n$ following lemma \ref{X47}, 
so if we set $\mu=\mp e^{-i{{\tau_n}\over 2}}P^\gamma_1\mu'$ it ensures 
$\mu \in E^n_{12}$ if $\mu' \in E^n_{21}$. Moreover, by definition of 
$\mu'$ we have: $$\mp e^{-i{{\tau_n}\over 2}}P^\gamma_2\mu=
 e^{-i{{\tau_n}}}P^\gamma_2 P^\gamma_1\mu'=\mu'
$$

\section {Domain Decomposition algorithm for  the Helmholtz equation}
\medskip

\subsection{The Domain Decomposition framework for  the Helm\-holtz equation}
\medskip

\bpr \label{X51}:  Let $k\in R$, $k\neq 0$. For any  $f  \in L^2(\Omega)$, let $u \in H^1_0(\Omega)$ be the unique solution of the Helmholtz equation 
$$\Delta u + k^2 u=f  $$
Let $\gamma \in R, \gamma \neq 0$. For $m=1,2$, let $v_m  \in H_m$ solve the equation
$$\Delta v_m + k^2 v_m=f_{\vert \Omega_m}\quad {\rm and}   \quad {\partial v_m \over \partial n_m}-i\gamma j \rho^\Gamma v_m =0 \,\, {\rm on}\,\,\Gamma$$
Let  $$\nu_m= 2i\gamma j \rho^\Gamma v_m \quad {\rm and} \quad \eta=(P_1^\gamma \nu_2, P_2^\gamma \nu_1)$$ 
Then the equation in $L^2(\Gamma)\times L^2(\Gamma)$ $$(A^\gamma -Id)\pi = \eta$$
has a unique solution: $$\pi=(\pi_1, \pi_2)\quad {\rm with} \quad\pi_m=  {\partial u \over \partial n_m}+i\gamma j \rho^\Gamma u-2i\gamma j \rho^\Gamma v_m$$
\epr
\medskip

\noindent {\it proof:} First notice that the assumption $f\in L^2(\Omega)$ and assumption {\bf (A)} on $k$  imply $u\in L^2(\Omega)$, so $\Delta u\in L^2(\Omega)$. Regularity of $\partial \Omega$ enables the use of classical regularity results ({\bf [Ag]}) for solutions of elliptic boundary problems to have $u\in H^2(\Omega)$, hence  ${\partial u \over \partial n_m }\in H^{1\over 2}(\Gamma) \subset L^2(\Gamma)$.  This proves that $\pi_m\in L^2(\Gamma)$. 

\noindent Let $w_m=u_{\vert \Omega_m}-v_m$. Then  $$\pi_m={\partial u \over \partial n_m}+i\gamma j \rho^\Gamma u-2i\gamma j \rho^\Gamma v_m= {\partial w_m \over \partial n_m}+i\gamma j \rho^\Gamma w_m$$ Because the function $w_m$ fulfills
$$\Delta w_m + k^2 w_m=0,  \quad w_m \in H_m, \quad {\partial w_m \over \partial n_m}-i\gamma j \rho^\Gamma w_m = {\partial u \over \partial n_m}-i\gamma j \rho^\Gamma u \,\,\, {\rm on}\,\,\, \Gamma$$
one has
$$P_m^\gamma ({\partial u \over \partial n_m}-i\gamma j \rho^\Gamma u)=\pi_m$$
so, if $m'=2,1$ for $m=1,2$
$$\pi_m=-P_m^\gamma ({\partial u \over \partial n_{m'}}+i\gamma j \rho^\Gamma u)=-P_m^\gamma (\pi_{m'}+\nu_{m'})$$
This writes
$$\pi=A^\gamma \pi - \eta$$
Uniqueness follows from theorem \ref{X50} (iii).
\medskip

\bpr \label{X52}:  Let $k\in R$, $k\neq 0$, $\gamma \in R, \gamma \neq 0$.  For  $f\in L^2(\Omega)$ and for $m=1,2$, let $v_m\in H_m$ solve the equation
$$\Delta v_m + k^2 v_m=f_{\vert \Omega_m}  \quad {\rm and} \quad {\partial v_m \over \partial n_m}-i\gamma j \rho^\Gamma v_m =0 \,\,\, {\rm on}\,\,\, \Gamma$$
Let  $$\nu_m= 2i\gamma j \rho^\Gamma v_m \quad {\rm and} \quad \eta=(P_1^\gamma \nu_2, P_2^\gamma \nu_1)$$ 
Let $\pi=(\pi_1, \pi_2)\in L^2(\Gamma)\times L^2(\Gamma)$ solve the equation:
$$(A^\gamma-Id) \pi = \eta$$
and let $u_m \in H_m$ solve the equation
$$\Delta u_m + k^2 u_m=f_{\vert \Omega_m} \quad {\rm and}\quad {\partial u_m \over \partial n_m}+i\gamma j \rho^\Gamma u_m =\pi_m+\nu_m \,\,\,{\rm on}\,\,\, \Gamma$$
Then $u$ given by $u_{\vert \Omega_m}=u_m$ solve the Helmholtz equation
$$\Delta u + k^2 u=f \quad u \in H^1_0(\Omega)$$
\epr
\medskip

\noindent {\it proof:}  By definition of $u_m$ and $v_m$ one has:
$$\Delta (u_m-v_m) + k^2 (u_m-v_m)=0,  \quad u_m-v_m \in H_m$$
and
$$ {\partial (u_m-v_m) \over \partial n_m}+i\gamma j \rho^\Gamma (u_m-v_m) =\pi_m \quad {\rm on}\quad \Gamma$$
This implies through proposition \ref{X30}:
$$(P_m^\gamma )^{-1}\pi_m=P_m^{-\gamma }\pi_m= {\partial (u_m-v_m) \over \partial n_m}-i\gamma j \rho^\Gamma (u_m-v_m) = {\partial u_m \over \partial n_m}-i\gamma j \rho^\Gamma u_m
$$
Because $\pi=A^\gamma \pi - \eta$ this implies (with  $m'=2,1$ for $m=1,2$)
$$ {\partial u_m \over \partial n_m}-i\gamma j \rho^\Gamma u_m = -\pi_{m'}-\nu_{m'}=
-{\partial u_{m'} \over \partial n_{m'}}-i\gamma j \rho^\Gamma u_{m'}$$ 
Adding and substracting these equalities gives:
$${\partial u_m \over \partial n_m}=
- {\partial u_{m'} \over \partial n_{m'}}\quad {\rm and}\quad \rho^\Gamma u_m=
\rho^\Gamma u_{m'}
$$
These jump conditions through $\Gamma$ imply that $\Delta u + k^2 u=f \,\, {\rm on}\,\, \Omega$,
and $u$ fulfills the Dirichlet boundary condition on $\partial \Omega$ because $u_m\in H_m$.

\brm \label{X53} Theorem \ref{X50} shows that the problem $$(\tilde A-Id)\pi =\eta
$$
is ill-posed for $\eta \in L^2(\Gamma)\times L^2(\Gamma)$. Proposition  \ref{X51} shows that if $\eta$ has the specific form given through the domain decomposition setting for  the Helmholtz equation, the equation $(\tilde A-Id)\pi =\eta$ do have a solution, (and this solution is unique by Theorem \ref{X50}). Proposition  \ref{X52} shows that this solution provides the solution of the Helmholtz equation.
\erm

\subsection {The domain decomposition $\theta$-algorithm for Helmholtz equation}

\noindent Let $f \in L^2(\Omega)$. Let $u\in H^1_0$ fulfill the non-dissipating Helmholtz equation $\Delta u + k^2 u=f$ in $\Omega$. The classical algorithm used, (for dissipating cavities with Sommerfeld-like boundary condition),  to solve by a domain decomposition technique the Helmholtz equation ({\bf [B],[BD],[D1],[D2],[CGJ]}) writes, in the non-dissipating case that discussed here, as follows:
for any $\pi^0=(\pi_1^0,\pi_2^0)$ given in $L^2(\Gamma)\times L^2(\Gamma)$ let $$\pi^{p+1}=\theta \pi^{p}+(1-\theta) A^\gamma \pi^{p}-(1-\theta)\eta
$$
where $\eta=(P_1^\gamma \nu_2, P_2^\gamma \nu_1)$ with 
$\nu_m= 2i\gamma j \rho^\Gamma v_m$ for $v_m \in H_m$ solving 
$$\Delta v_m + k^2 v_m=f_{\vert \Omega_m}  \quad {\rm and} \quad {\partial v_m \over \partial n_m}-i\gamma j \rho^\Gamma v_m =0 \,\,\, {\rm on}\,\,\, \Gamma$$
It is straightforward  to translate the $\theta$-algorithm in a PDE setting: one  use theorem \ref{X28} to get  the (unique) function $w^n_m \in H_m$ such that 
$$ \Delta w_m^{p} +k^2 w_m^{p}=0\quad {\rm and}\quad \pi^{p}_{m}= {\partial w_m^{p} \over \partial n_m}+i\gamma j\rho^\Gamma w_m^{p}
$$
In these $w^p=(w^p_1,w^p_2)$ variables the  $\theta$-algorithm  becomes:
 for $m=1$ and $2$ (resp. $m'=2$ and $1$)
$$ \Delta w_m^{p+1} +k^2 w_m^{p+1}=0, \quad w_m^{p+1}\in H_m $$
$$ {\partial w_m^{p+1} \over \partial n_m}-i\gamma j\rho^\Gamma w_m^{p+1}=\theta [ {\partial w_m^{p} \over \partial n_m}-i\gamma j\rho^\Gamma w_m^{p}]-(1-\theta)[{\partial w_{m'}^{p} \over \partial p_{m'}}+i\gamma j\rho^\Gamma w_{m'}^{p}+\nu_{m'}]\,\, {\rm on}\,\, \Gamma
$$
\medskip

\noindent For practical use in computing codes, one writes this algorithm in the $u^n=(u^n_1,u^n_2)$ variables with $u^p_m=w^p_m+v_m$  and gets: 
$$ \Delta u_m^{p+1} +k^2 u_m^{p+1}=f_{\vert \Omega_m}, \quad u_m^{p+1}\in H_m $$
$$ {\partial u_m^{p+1} \over \partial n_m}-i\gamma j\rho^\Gamma u_m^{p+1}=\theta [ {\partial u_m^{p} \over \partial n_m}-i\gamma j\rho^\Gamma u_m^{p}]-(1-\theta)[ {\partial u_{m'}^{p} \over \partial n_{m'}}+i\gamma j\rho^\Gamma u_{m'}^{p}]\,\, {\rm on}\,\, \Gamma
$$
\subsection{Convergence results for the $\theta$-algorithm}
Notice that 
if the sequence $(\pi^p)_{p\in N}$ has a limit $\pi^\infty$ in $L^2(\Gamma)\times L^2(\Gamma)$  then continuity of $A^\gamma$ gives:
$$(A^\gamma -Id)\pi^\infty=\eta
$$
and $\pi^\infty$ provides the solution of the Helmholtz equation on $\Omega$ as stated in proposition \ref{X52}.

\noindent Alternatively, a way to solve the Helmholtz equation on $\Omega$ through sol\-ving Helmholtz equations on $\Omega_m$, ($m=1,2$),  is to notice that  convergence of $(\pi^p)_{p\in N}$ in $L^2(\Gamma)\times L^2(\Gamma)$ implies convergence in $\Lambda'\times \Lambda'$, and theorem \ref{X28} shows that the sequence $(w^p)_{{p\in N}}=(\,(w^p_1,w^p_2)\,)_{{p\in N}}$ has a limit in $H_1\times H_2$, which implies convergence of the sequence $(u^p)_{{p\in N}}=(\,(u^p_1,u^p_2)\,)_{{p\in N}}$ in $H_1\times H_2$. Proposition \ref{X52} shows that its limit $u^\infty=(u^\infty_1,u^\infty_2)$ provides the solution of the Helmholtz equation, through.
$$(u^\infty_1,u^\infty_2)=(u_{\vert \Omega_1},u_{\vert \Omega_2})
$$
Here are the convergence results for the $\theta$-algorithm. We begin with a negative result:

\bpr \label{X55} If $\theta =0$ then  the sequence  $(\pi^p)_{p\in N}$ has no limit in $L^2(\Gamma)\times L^2(\Gamma)$ unless if its initial value fulfills $(A^\gamma -Id)\pi^0=\eta$. Written for the $(u^p)_{p\in N}$ sequence, this translates to $u^0_1=u_{\vert \Omega_1}$ and 
$u^0_2=u_{\vert \Omega_2}$ where $u\in H^1_0(\Omega)$ solves $\Delta u + k^2 u=f$ in $\Omega$
  \epr
\noindent {\it proof:} for $\theta =0$ one has $\pi^p-\pi^{p-1}=A^\gamma(\pi^{p-1}-\pi^{p-2})$ and proposition \ref{X35} gives
$$\forall p \quad\Vert \pi^p-\pi^{p-1} \Vert_{L^2(\Gamma)\times L^2(\Gamma)}=\Vert \pi^{p-1}-\pi^{p-2} \Vert_{L^2(\Gamma)\times L^2(\Gamma)}
$$
This prevents convergence unless if $\pi^1=\pi^{0}$, i.e. $\pi^{0}$ fulfills 
$$\pi^{0}= A^\gamma \pi^{0}-\eta
$$
i.e. unless $u^0_1=u_{\vert \Omega_1}$ and 
$u^0_2=u_{\vert \Omega_2}$ by proposition \ref{X52}.

\brm \label{X56}  If $\theta =0$  the sequence $(u^p)_{{p\in N}}$ may have a limit in $H_1\times H_2$ even if $u^0_1\neq u_{\vert \Omega_1}$ or 
$u^0_2\neq u_{\vert \Omega_2}$. This is because convergence of $(u^p)_{{p\in N}}$ in $H_1\times H_2$ implies convergence of 
$( {\partial u_{1}^{p} \over \partial n_{1}}+i\gamma j\rho^\Gamma u_{1}^{p} ,{\partial u_{2}^{p} \over \partial n_{2}}+i\gamma j\rho^\Gamma u_{2}^{p})$ in $\Lambda' \times  \Lambda'$, i.e. convergence of $( {\partial w_{1}^{p} \over \partial n_{1}}+i\gamma j\rho^\Gamma u_{1}^{p} ,{\partial w_{2}^{p} \over \partial n_{2}}+i\gamma j\rho^\Gamma u_{2}^{p})$ in $\Lambda' \times  \Lambda'$, i.e. convergence of $(\pi^p)_{p\in N}$  in $\Lambda' \times  \Lambda'$, which do not contradict divergence  in  $L^2(\Gamma)\times L^2(\Gamma)$.
\erm
\medskip

\noindent We now turn to the main result:

\bth :\label{X57} For $f \in L^2(\Omega)$ and $k\in R$ let $u \in H^1_0(\Omega)$ solve $\Delta u+k^2u=f$. Let $\gamma \neq 0, \gamma\in R$.
Then for any $0<\theta <1$:
 \medskip
 
\noindent (i)  the sequence
$(\pi^p)_{p\in N}$ given by the
$\theta$-algorithm  converge in $L^ 
2(\Gamma) \times L^ 
2(\Gamma)$ to  $\pi^u=(\pi_1^u,\pi_2^u)$ with:
$$\pi_1^u={{\partial u} \over {\partial n_1}}+i\gamma j \rho^\Gamma u
\quad {\rm and}\quad \pi_2^u={{\partial u} \over {\partial n_2}}+i\gamma j \rho^\Gamma u$$

\noindent (ii)  the sequence
$(u^p)_{p\in N}$ given by the
$\theta$-algorithm converge in $H_1\times H_2$ to $(u_{\vert
\Omega_1},u_{\vert \Omega_2})$ 
\medskip

\noindent (iii) There is no uniform geometric rate of convergence.
\eth
\medskip

\noindent {\it Proof:} 

\noindent (i) let $\pi^u=(\pi_1^u,\pi_2^u)$ with:
$$\pi_1^u={{\partial u} \over {\partial n_1}}+i\gamma j \rho^\Gamma u
\quad {\rm and}\quad \pi_2^u={{\partial u} \over {\partial n_2}}+i\gamma j \rho^\Gamma u$$
Nullity of jumps of $u$ and its normal derivatives  through $\Gamma$  gives
$$\pi^u =A\pi^u+\eta$$
This implies $$ \pi^p-\pi^u=\theta [\pi^{p-1}-\pi^u]+(1-\theta)A[\pi^{p-1}-\pi^u] $$
We  use   eigenprojectors of $A^\gamma$ given by theorem \ref{X50} and denote by:
$$\delta^p_{n,\pm}=P_n^\pm (\pi^p-\pi^u)$$
Completeness of the set of orthogonal eigenprojectors $(\,P_n^\pm\,)_{n,\pm}$  proved in theorem \ref{X50} gives:
$$ \pi^p-\pi^u=\sum_n \delta^p_{n,+}+\sum_n \delta^p_{n,-}
$$
Decomposition of $L^2(\Gamma) \times L^2(\Gamma) $ by eigenspaces of $A^\gamma$ writes for successive terms of the $\theta-$algorithm sequence as follows:
$$ \delta_{n,\pm}^{p}=[\theta \pm (1-\theta)e^{i{\tau_n \over 2}}] \delta_{n,\pm}^{p-1}$$
This implies:
$$\Vert \delta_{n,\pm}^{p}\Vert_{ L^2(\gamma) \times L^2(\gamma)}=[1-2\theta (1-\theta)(1 \mp \cos{\tau_n \over 2})]^p\Vert  \delta_{n,\pm}^{0}\Vert_{ L^2(\gamma) \times L^2(\gamma)}
$$ 
and orthogonality of the eigenprojectors writes:
$$\Vert \pi^p-\pi^u \Vert^2_{ L^2(\gamma) \times L^2(\gamma)}=
\sum_{n,\pm} [1-2\theta (1-\theta)(1 \mp \cos{\tau_n \over 2})]^{2p} \Vert \delta_{n,\pm}^{0} \Vert^2_{ L^2(\gamma) \times L^2(\gamma)}
$$
Assumption $0<\theta <1$ and theorem \ref{X50} (with proposition \ref{X43} asserting $\tau_n \neq 0$ mod $2\pi$) imply $$0<1-2\theta (1-\theta)(1 \mp \cos{\tau_n \over 2})<1$$ 
and Lebesgue convergence theorem gives
$$\pi^p \,\, \stackrel{{ L^2(\Gamma) \times L^2(\Gamma)}}{\longrightarrow}\,\, \pi^u 
$$

\noindent This proves assertion (i).
\medskip

\noindent Assertion (ii) is straightforward: convergence of $\pi^p$ to $\pi^u$ in  $L^ 
2(\Gamma) \times L^ 
2(\Gamma)$ implies its convergence in $\Lambda' \times \Lambda'$ which implies convergence of the related sequence $w^p$ in $H_1\times H_2$, and accordingly convergence of $u^p$ to $(u_{\vert \Omega_1},u_{\vert \Omega_2})$
\medskip

\noindent  Assertion (iii) is obvious by taking initial data for $\pi^p$ in the $n-$th eigenspace of $A^\gamma$ and notice that $\tau_n \longrightarrow 0$ (proposition \ref{X44})

\section {Bibliography}
\medskip

\noindent {\bf [AG] Agmon S..} Lectures on elliptic boundary value problems. Van Nostrand - Math Studies n¡2 - 1965
\medskip

\noindent {\bf [B] Benamou J.D.} A domain decomposition method for the optimal control of system governed by the Helmholtz equation. {\it 3rd international conference on math. and num. aspects of  waves propagation phenomena. SIAM ed} - Cannes-Mandelieu -(1995)
\medskip

\noindent {\bf [BD] Benamou J.D. -Despres B.} A domain decomposition method for the  Helmholtz equation and related optimal control problems. {\it J. of Comp. Physics} vol 136  (1995) pp 68-82
\medskip

\noindent {\bf [CGJ] Collino F. - Ghanemi S. - Joly P.} Domain decomposition method for harmonic wave propagation: a general presentation {\it Rapport de recherche INRIA 3473} (1998)
\medskip

\noindent {\bf [D1] Despres B.} Domain decomposition method and the  Helmholtz problem {\it Math. and num. aspects of waves propagation phenomena.  SIAM ed} (1993) 
\medskip

\noindent {\bf [D2] Despres B.} Domain decomposition method and the  Helmholtz problem (part ii) {\it 2nd international conference on math. and num. aspects of waves propagation phenomena. SIAM ed}  (1995)
\medskip

\noindent {\bf [FMQT] Faccioli E. - Maggio F. - Quarteroni A. - Tagliani A.} Spectral domain decomposition methods for the solution of elastice and acoustic wave equation. {\it Geophysics 61}  (1996) pp 1160-1174
\medskip

\noindent {\bf [FQZ] Funaro D. - Quarteroni A. - Zanolli P.} An iterative procedure with interface relaxation for domain decomposition methods {\it SIAM J. Numer. Anal.} 25 (1988) pp 1213-1236
\medskip

\noindent {\bf [L] Lions P.L.} On the Schwarz alternating method I. {\it First International Symposium on domain decomposition methods for PDE (1988) R.Glowinski et als eds. SIAM.} Philadelphia. pp 1-42
\medskip

\noindent {\bf [PW] Proskwowski W. - Widlund O.} On the numerical solutions of Helmholtz's equation by the capacitance matrix method {\it Math. Comp.} (1976) vol 30 pp 433-468

\medskip
\noindent {\bf [QV] Quarteroni A. - Valli A.} Domain decomposition methods for PDE. {\it Num. Math. and  Sc. Comp. Clarendon Press.} Oxford (1999)
\medskip

\end{document}